\documentclass[11pt]{article}
\usepackage{amssymb,latexsym}
\usepackage{epsfig}
\usepackage{eufrak}
\usepackage{amsmath}
\usepackage{mathrsfs}
\usepackage{color}

\newtheorem{theorem}{Theorem}
\newtheorem{lemma}{Lemma}
\newtheorem{corollary}{Corollary}

\newcommand{\be}{\begin{equation}}
\newcommand{\ee}{\end{equation}}
\newcommand{\bee}{\begin{eqnarray*}}
\newcommand{\eee}{\end{eqnarray*}}
\newcommand{\bel}{\begin{eqnarray}}
\newcommand{\eel}{\end{eqnarray}}
\newcommand{\bec}{\begin{cases}}
\newcommand{\eec}{\end{cases}}
\newcommand{\bem}{\begin{bmatrix}}
\newcommand{\eem}{\end{bmatrix}}

\newcommand{\la}{\label}
\newcommand{\li}{\left}
\newcommand{\ri}{\right}

\newcommand{\DEF}{\stackrel{\mathrm{def}}{=}}

\newcommand{\ovl}{\overline}

\newcommand{\vep}{\varepsilon}
\newcommand{\lm}{\lambda}
\newcommand{\Lm}{\Lambda}

\newcommand{\si}{\sigma}

\newcommand{\de}{\delta}

\newcommand{\vDe}{\varDelta}

\newcommand{\ga}{\gamma}
\newcommand{\Ga}{\Gamma}
\newcommand{\vse}{\vartheta}
\newcommand{\se}{\theta}
\newcommand{\Se}{\Theta}

\newcommand{\ze}{\zeta}
\newcommand{\al}{\alpha}

\newcommand{\ro}{\rho}

\newcommand{\Om}{\Omega}

\newcommand{\f}{\frac}
\newcommand{\sq}{\sqrt}
\newcommand{\cd}{\cdots}

\newcommand{\qu}{\quad}
\newcommand{\qqu}{\qquad}

\newcommand{\mscr}{\mathscr}
\newcommand{\mcal}{\mathcal}
\newcommand{\mbf}{\mathbf}
\newcommand{\bb}{\mathbb}

\newcommand{\wh}{\widehat}

\newcommand{\mrm}{\mathrm}
\newcommand{\bs}{\boldsymbol}

\newcommand{\sh}{\slash}

\newcommand{\tx}{\text}

\newcommand{\iy}{\infty}

\newcommand{\bed}{\begin{description}}
\newcommand{\eed}{\end{description}}
\newcommand{\bei}{\begin{itemize}}
\newcommand{\eei}{\end{itemize}}
\newcommand{\ben}{\begin{enumerate}}
\newcommand{\een}{\end{enumerate}}
\newcommand{\bib}{\bibitem}
\newcommand{\beL}{\begin{lemma}}
\newcommand{\eeL}{\end{lemma}}
\newcommand{\beT}{\begin{theorem}}
\newcommand{\eeT}{\end{theorem}}
\newcommand{\sect}{\section}

\newcommand{\bpf}{\begin{pf}}
\newcommand{\epf}{\end{pf}}

\newcommand{\bi}{\binom}

\newcommand{\pfbox}{\hfill\mbox{$\Box$}}

\newenvironment{pf}{\paragraph*{Proof{\rm.}}}{\pfbox\bigskip}

\begin{document}

\title{{\bf New Probabilistic Inequalities from Monotone Likelihood Ratio Property}
\thanks{The author had been previously working with Louisiana
State University at Baton Rouge, LA 70803, USA, and is now with
Department of Electrical Engineering, Southern University and A\&M
College, Baton Rouge, LA 70813, USA; Email: chenxinjia@gmail.com}}

\author{Xinjia Chen}

\date{First submitted in June 2010}

\maketitle

\begin{abstract}

In this paper, we propose a new approach for deriving probabilistic
inequalities.  Our main idea is to exploit the information of
underlying distributions by virtue of the monotone likelihood ratio
property and Berry-Essen inequality. Unprecedentedly sharp bounds
for the tail probabilities of some common distributions are
established.  The applications of the probabilistic inequalities in
parameter estimation are discussed.

\end{abstract}

\section{Introduction}  \la{Intr}

Probabilistic inequalities are important ingredients of fundamental
probabilisty theory.  A classical approach for deriving
probabilistic inequalities is based on the moment or moment
generating functions of relevant random variables.  In view of the
fact that the moment generating function is actually a  moment
function in a general sense, we call this approach as {\it Method of
Moments}. Many well-known inequalities such as Markov inequality,
Chebyshev inequality, Chernoff bounds \cite{Chernoff}, Hoeffding
\cite{Hoeffding} inequalities are developed in this framework. In
order to use the method of moments to derive probabilistic
inequalities, a critical step is to obtain a closed-form expression
for the moment or moment generating function. However, for some
common distributions, the moment or moment generating function may
be either unavailable or too complicated for analytical treatment.
Familiar examples are Student's $t$-distribution, Snedecor's
$F$-distribution, hypergeometric distribution, hypergeometric
waiting-time distribution, for which the method of moments is not
useful for deriving sharp bounds for tail probabilities. In addition
to this limitation, another drawback of the method of moments is
that the information of the underlying distribution may not be fully
exploited.  This is especially true when the relevant distribution
is analytical and known.

In this paper, we take a new path to derive probabilistic
inequalities.  In order to overcome the limitations of the method of
moments, we exploit the information of underlying distribution by
virtue of the statistical concept of Monotone Likelihood Ratio
Property (MLRP).  We discovered that, the MLRP is extremely powerful
for deriving sharp bounds for the tail probabilities of a large
class of distributions. Specially, in combination of the Berry-Essen
inequality, the MLRP can be employed to improve upon the
Chernoff-Hoeffding bounds for the tail probabilities of the
exponential family by a factor about two. For common distributions
such as Student's $t$-distribution, Snedecor's $F$-distribution,
hypergeometric distribution, hypergeometric waiting-time
distribution, we also obtained unprecedentedly sharp bounds for the
tail probabilities.  We demonstrate that the MRLP can be used to
illuminate probabilistic phenomenons with very elementary knowledge.

The remainder of the paper is organized as follows.  In Section 2, we present our most general results, especially the Likelihood Ratio Bounds
(LRB).  Section 3 gives bounds on the distribution of likelihood ratio. In Section 4, we develop a unified theory for bounding the tail
probabilities of the exponential family. In Section 5, we apply our general theory to obtain tight bounds for the tail probabilities of common
distributions. In Section 6, we explore the general applications of the probabilistic inequalities for parameter estimation.  Section 7 is the
conclusion.  Throughout this paper, we shall use the following notations.  The set of real numbers is denoted by $\bb{R}$. The set of integers
is denoted by $\bb{Z}$.   We use the notation $\Pr \{ . \mid \se \}$ to indicate that the associated random samples $X_1, X_2, \cd$ are
parameterized by $\se$. The parameter $\se$ in $\Pr \{ . \mid \se \}$  may be dropped whenever this can be done without introducing confusion.
The expectation of a random variable is denoted by $\bb{E}[.]$.  The notation $I_Z$ denotes the support of $Z$. The other notations will be made
clear as we proceed.

\sect{Likelihood Ratio Bounds} \la{secLRB}

The statistical concept of monotone likelihood ratio plays a central role in our development of new probabilistic inequalities. Before
presenting our new results, we shall describe the MLRP as follows. Let $X_1, X_2, \cd, X_n$ be a sequence of random variables defined in
probability space $(\Om, \mscr{F}, \Pr )$ such that the joint distribution of $X_1, \cd, X_n$ is determined by parameter $\se$ in $\Se$.  Let
$f_n(x_1, \cd, x_n; \se)$ be the joint probability density function for the continuous case or the probability mass function for the discrete
case, where $(x_1, \cd, x_n)$ denotes a realization of $(X_1, \cd, X_n)$.  The family of joint probability density or mass functions is said to
posses MLRP if there exist a nonnegative multivariate function $\Lm(z, \vse_0, \vse_1)$ of $z \in \mscr{Z}, \; \vse_0 \in \Se, \; \vse_1 \in
\Se$ and a multivariate function $\varphi = \varphi (x_1, \cd, x_n)$ of $x_1, \cd, x_n$ such that the following requirements are satisfied.

(I) $\varphi = \varphi (x_1, \cd, x_n)$ takes values in $\mscr{Z}$ for arbitrary realization, $(x_1, \cd, x_n)$, of $(X_1, \cd, X_n)$.

(II) For arbitrary parametric values $\se_0,  \se_1 \in \Se$, the function $\Lambda (z, \se_0, \se_1)$ is non-decreasing with respect to $z \in
\mscr{Z}$ provided that $\se_0 \leq \se_1$.

(III)  For arbitrary parametric values $\se_0, \se_1 \in \Se$, the likelihood ratio $\f{ f_n  (x_1, \cd, x_n; \se_1) } { f_n (x_1, \cd, x_n;
\se_0)  }$ can be expressed as $\Lambda (\varphi, \se_0, \se_1)$.

Now we are ready to state our general results as Theorem \ref{THm1}
in the following.

\beT   \la{THm1} Let $\bs{\varphi} = \varphi(X_1, \cd, X_n)$.  Let $\vse (z)$ be a function of $z \in \mscr{Z}$ taking values in $\Se$. Suppose
the monotone likelihood ratio property holds. Define $\mscr{M} (z, \se ) = \Lambda (z, \vse(z), \se)$ for $z \in \mscr{Z}$ and $\se \in \Se$.
Then, \be \la{bound1}
 \Pr \{ \bs{\varphi} \geq z \mid \se \} \leq \mscr{M} (z, \se ) \times
  \Pr \{ \bs{\varphi} \geq z \mid \vse (z)  \} \leq \mscr{M} (z, \se ) \ee
  for $z \in \mscr{Z}$ such that $\vse (z)$ is no less than $\se \in
\Se$. Similarly, \be \la{bound2}
 \Pr \{ \bs{\varphi} \leq z \mid \se \} \leq \mscr{M} (z, \se )
 \times \Pr \{ \bs{\varphi} \leq z \mid \vse (z) \} \leq \mscr{M} (z, \se )  \ee
 for $z \in \mscr{Z}$ such that $\vse (z)$ is no greater than $\se
\in \Se$.

Assume that the following additional assumptions are satisfied:

(a) $\vse (z) = z$ for any $z \in \Se$;

(b) $f_n (x_1, \cd, x_n; \se)$ can be expressed as a function $g(\varphi, \se)$ of $\varphi = \varphi(x_1, \cd, x_n)$ and $\se$;

(c) $g(z, \se)$ is non-decreasing with respect to $\se \in \Se$ no greater than $z \in \Se$ and is non-increasing with respect to $\se \in \Se$
no less than $z \in \Se$.

Then,  the following statements hold true:

 (i) $\mscr{M} (z, \se ) = \Lambda (z, z, \se) = \f{ g(z, \se)  }
 { g(z, z) }$ for $z, \se \in \Se$.

 (ii) $\mscr{M} (z, \se )$ is non-decreasing with respect to $\se \in
 \Se$ no greater than $z \in \Se$ and is non-increasing with respect to $\se \in
 \Se$ no less than $z \in \Se$.

 (iii) $\mscr{M} (z, \se )$ is non-decreasing with respect to $z \in \Se$
 no greater than $\se \in \Se$ and is non-increasing with respect
to $z \in \Se$ no less than $\se \in \Se$.

 \eeT

 The  proof of Theorem \ref{THm1} is provided in Appendix \ref{THm1_app}.
 Since inequalities (\ref{bound1}) and (\ref{bound2}) are derived
 from the MLRP, these inequalities are referred to as the {\it Likelihood
 Ratio Bounds} in this paper and its previous version \cite{Chen2}.

 An immediate application of Theorem \ref{THm1} can be found in the area of statistical hypothesis testing.  It is a frequent problem to test hypothesis $\mscr{H}_0: \se \leq \se_0$ versus $\mscr{H}_1: \se \geq \se_1$,
 where $\se_0 < \se_1$ are two parametric values in $\Se$. Assume that there is a statistic $\wh{\se}$ defined in terms of $X_1, \cd, X_n$ such
 that the probability ratio $\f{ f_n (X_1, \cd, X_n; \se_1)  }{ f_n (X_1, \cd, X_n; \se_0) }$ can be expressed as $\Lm (\wh{\se}, \se_0,
 \se_1)$, which is increasing with respect to $\wh{\se}$.  To test the hypotheses, a classical method is to choose a number $\ga \in \Se$ such that
$\se_0 \leq \ga \leq \se_1$ and make the decision: Accept $\mscr{H}_0$ if $\wh{\se} \leq \ga$ and otherwise reject $\mscr{H}_0$.  To offer
simple bounds for the risks of making an erroneous decision, we have obtained the following new result: \be \la{testbound} \Pr \{ \tx{Reject} \;
\mscr{H}_0 \mid \mscr{H}_0 \} \leq \Lm (\ga, \ga, \se_0), \qqu \Pr \{ \tx{Reject} \;  \mscr{H}_1 \mid \mscr{H}_1 \} \leq \Lm (\ga, \ga,
 \se_1).
\ee To prove (\ref{testbound}), note that
\[
\Pr \{ \tx{Reject} \;  \mscr{H}_0 \mid \mscr{H}_0 \} = \Pr \{ \wh{\se} > \ga \mid \mscr{H}_0 \} \leq \Pr \{ \wh{\se} > \ga \mid \se_0 \} \leq
\Lm (\ga, \ga, \se_0),
\]
where the first inequality is due to the monotonicity of the likelihood ratio, and the second inequality is a consequence of Theorem \ref{THm1}.
Similarly, \[ \Pr \{ \tx{Reject} \;  \mscr{H}_1 \mid \mscr{H}_1 \} = \Pr \{ \wh{\se} \leq \ga \mid \mscr{H}_1 \} \leq \Pr \{ \wh{\se} \leq \ga
\mid \se_1 \} \leq \Lm (\ga, \ga, \se_1).
\]
It can be checked that such bounds apply to the exponential family and hypergeometric distribution.

\section{Bounds on the Distribution of Likelihood Ratio} \la{secCDF}

Let $f_X(x; \se)$ denote the probability density (or mass) function of $X$ parameterized by $\se \in \Se$.   Let $X_1, X_2, \cd$ be i.i.d.
samples of $X$.  Consider hypothesis $\mscr{H} : \se = \se_0$. Assume that for a sample of size $n$, there exists a maximum likelihood estimator
(MLE) $\wh{\se}_n$ for $\se_0$ such that the sequence of estimators $\wh{\se}_n, \; n = 1, 2, \cd$ converges in probability to $\se_0$. Define
likelihood ratio
\[
\lm_{\mscr{H}} = \f{ \prod_{i=1}^n f_X(X_i; \se_0)  }{ \prod_{i=1}^n f_X(X_i; \wh{\se}_n)  }, \qqu n = 1, 2, \cd.
\]
Assume that $\wh{\se}_n$ is asymptotically normally distributed with mean $\se_0$. In this setting, Wilks proved that
\[
\lim_{n \to \iy} \Pr \{  - 2 \ln  \lm_{\mscr{H}} < \chi^2 \mid \se_0 \} = \f{1}{\sq{2 \pi}} \int_0^{\chi^2} u^{-\f{1}{2}} e^{-\f{u}{2} } d u
\]
that is, if $\mscr{H}$ is true, $- 2 \ln  \lm_{\mscr{H}}, \; n = 1, 2, \cd$ converges in distribution to the chi-square distribution of degree
one.   The proof of this result can be found in pages 410--411 of Wilks' text book {\it Mathematical Statistics}.  This result has important
application for testing hypothesis $\mscr{H} : \se = \se_0$.  Suppose the decision rule is that $\mscr{H}$ is rejected if $- 2 \ln
\lm_{\mscr{H}} > \chi_\al^2$, where $\chi_\al^2$ is the number for which $\Pr \{ \chi^2 > \chi_\al^2 \} = \al$.  Then, $\lim_{n \to \iy} \Pr \{
\tx{Reject} \; \mscr{H} \mid  \mscr{H} \} = \al$.

The drawback of the asymptotic result is that it is not clear how large the sample size $n$ is sufficient for the asymptotic distribution to be
applicable. To address this issue, it is desirable to obtain tight bounds for the distribution of $- 2 \ln  \lm_{\mscr{H}}$. For this purpose,
we can apply Theorem \ref{THm1} to derive the following results.

\beT \la{boundCDFLR} Let $\al$ be a positive number and $n$ be a positive integer.  Let $f_n( x_1, \cd, x_n; \se)$ denote the joint probability
density or mass function of random variables $X_1, \cd, X_n$ parameterized by $\se \in \Se$. Assume that $\f{ f_n( X_1, \cd, X_n; \se_1) }{f_n(
X_1, \cd, X_n; \se_0)}$ can be expressed as a function, $\Lm (\varphi_n , \se_0, \se_1)$,  of $\se_0, \se_1$ and $\varphi_n = \varphi(X_1, \cd,
X_n)$ such that $\Lm (\varphi_n , \se_0, \se_1)$ is increasing with respect to $\varphi_n$. Let $\wh{\se}_n$ be a function of $\varphi_n$ such
that $\wh{\se}_n$ takes values in $\Se$. Then, \bel &  & \Pr \li \{ \f{ f_n( X_1, \cd, X_n; \se) }{f_n( X_1, \cd, X_n; \wh{\se}_n)} \leq
\f{\al}{2}, \; \wh{\se}_n \leq \se \mid \se \ri \}
\leq \f{\al}{2}, \la{RB1}\\
&  & \Pr \li \{ \f{ f_n( X_1, \cd, X_n; \se) }{f_n( X_1, \cd, X_n; \wh{\se}_n)} \leq \f{\al}{2}, \; \wh{\se}_n \geq \se \mid \se \ri \} \leq \f{\al}{2}, \la{RB2}\\
&  & \Pr \li \{ \f{ f_n( X_1, \cd, X_n; \se) }{f_n( X_1, \cd, X_n; \wh{\se}_n)} \leq \f{\al}{2} \mid \se \ri \} \leq \al \la{RB3} \eel for $\se
\in \Se$.  Moreover, under additional assumption that $\wh{\se}_n$ is a MLE for $\se$, the following inequalities \bel & & \Pr \li \{ \f{
\sup_{\vse \in \mscr{S}} f_n( X_1, \cd, X_n; \vse) }{\sup_{\vse \in \Se} f_n( X_1, \cd, X_n; \vse)} \leq \f{\al}{2}, \;
\wh{\se}_n \leq \inf \mscr{S} \mid \se \ri \} \leq \f{\al}{2} \la{LBA},\\
 &  & \Pr \li \{  \f{ \sup_{\vse \in \mscr{S}} f_n( X_1, \cd, X_n;
\vse) }{\sup_{\vse \in \Se} f_n( X_1, \cd,
X_n; \vse)} \leq \f{\al}{2}, \; \wh{\se}_n \geq \sup \mscr{S} \mid \se \ri \} \leq \f{\al}{2},  \la{LBB}\\
&  & \Pr \li \{  \f{ \sup_{\vse \in \mscr{S}} f_n( X_1, \cd, X_n; \vse) }{\sup_{\vse \in \Se} f_n( X_1, \cd, X_n; \vse)} \leq \f{\al}{2} \mid
\se \ri \} \leq \al \la{LBC}
 \eel hold true for arbitrary nonempty subset $\mscr{S}$ of $\Se$ and all
$\se \in \mscr{S}$.

\eeT

See Appendix \ref{boundCDFLR_app} for a proof.  To apply inequalities (\ref{RB1})--(\ref{RB3}), there is no necessity for  $X_1, \cd, X_n$ to be
i.i.d. and $\wh{\se}_n$ to be a MLE for $\se$.  Applying Theorem \ref{boundCDFLR} to the likelihood ratio
\[ \lm_{\mscr{H}} = \f{ f_n(X_1, \cd, X_n; \se_0)  }{  f_n(X_1, \cd, X_n; \wh{\se}_n)  }
\]
yields
\[
\Pr \{  - 2 \ln  \lm_{\mscr{H}} \geq  \chi^2 \mid \se_0 \} \leq 2 \exp \li (  - \f{\chi^2}{2} \ri ).
\]
As a by product, we have proved the inequality \[ \f{1}{\sq{2 \pi}} \int_{z}^\iy u^{-\f{1}{2}} e^{-\f{u}{2} } d u < 2 \exp \li ( - \f{z}{2} \ri
), \qqu z
> 0.
\]
With regard to testing hypothesis $\mscr{H} : \se = \se_0$, if the decision rule is to reject $\mscr{H}$ when $\lm_{\mscr{H}}  \leq \f{\al}{2}$,
then
\[
\Pr \{ \tx{Reject} \; \mscr{H} \mid  \mscr{H} \} = \Pr \li \{  \lm_{\mscr{H}}  \leq \f{\al}{2} \mid \se_0 \ri \} \leq \al.
\]
Since the acceptance region is
\[
\li \{ (x_1, \cd, x_n): \f{ f_n(x_1, \cd, x_n; \se_0)  }{  f_n(x_1, \cd, x_n; \wh{\se}_n)  } > \f{\al}{2} \ri \},
\]
it follows that inverting the acceptance region leads to a confidence region for $\se$ with coverage probability no less than $1 - \al$.
Specially, if we define random region \[ \mscr{R} = \li \{ \se_0 \in \Se: \f{ f_n(X_1, \cd, X_n; \se_0)  }{  f_n(X_1, \cd, X_n; \wh{\se}_n)  } >
\f{\al}{2} \ri \},
\]
then $\Pr \{ \se \in \mscr{R} \mid \se \} \geq 1 - \al$ for all $\se \in \Se$.  It can be shown that $\mscr{R}$ is actually an interval if
$\wh{\se}_n$ is a MLE for $\se$.  We will return to the problem of interval estimation later.

\sect{Probabilistic Inequalities for Exponential Family} \la{Prexp}

Our main objective for this section is to develop a unified theory
for bounding the tail probabilities of the exponential family.  A
single-parameter exponential family is a set of probability
distributions whose probability density function (or probability
mass function, for the case of a discrete distribution) can be
expressed in the form \be \la{expdef}
 f_X (x, \se) = h (x) \exp ( \eta (\se) T( x ) - A (\se) ), \qqu \se
 \in \Se
\ee where $T(x), \; h(x), \eta(\se)$, and $A(\se)$ are known
functions.

For the exponential family described above, we have the following
results.

\beT  \la{THm2} Let $X$ be a random variable with probability density function or probability mass function defined by (\ref{expdef}). Let $X_1,
\cd, X_n$ be i.i.d. samples of $X$. Define $\wh{\bs{\se}} = \f{ \sum_{i = 1}^n T(X_i) } { n }$ and $\mscr{M} (z, \se) = \li [ \f{ \exp ( \eta
(\se) z - A (\se) ) } { \exp ( \eta (z) z - A (z) ) } \ri ]^n$ for $z, \se \in \Se$. Suppose that $\f{d \eta (\se)}{d \se}$ is positive for $\se
\in \Se$. Then,
\[
\Pr \li \{ \wh{\bs{\se}} \geq z \mid \se \ri \} \leq  \mscr{M} (z, \se) \times \Pr \li \{ \wh{\bs{\se}} \geq z \mid z \ri \} \qqu \tx{for $z \in
\Se$ no less than $\se \in \Se$}
\]
and \[ \Pr \li \{ \wh{\bs{\se}} \leq z \mid \se \ri \} \leq \mscr{M} (z, \se) \times \Pr \li \{ \wh{\bs{\se}} \leq z \mid z \ri \}  \qqu \tx{for
$z \in \Se$ no greater than $\se \in \Se$}.
\]
Moreover, under the additional assumption that $\f{ d A ( \se ) }{ d
\se } = \se \f{ d \eta ( \se ) }{d \se}$, the following statements
hold true:

(i) $\wh{\bs{\se}}$ is a maximum-likelihood and unbiased estimator
of $\se$.

(ii) $\mscr{M} (z, \se) = \inf_{t \in \bb{R} } \bb{E} \li [ \exp \li
( n t ( \wh{\bs{\se}} - z) \ri ) \ri ]$, where the infimum is
attained at $t = \eta(z) - \eta (\se)$.

(iii) $\mscr{M} (z, \se )$ is increasing with respect to $\se \in \Se$ no greater than $z \in \Se$ and is decreasing with respect to $\se \in
\Se$ no less than $z \in \Se$.

(iv) $\mscr{M} (z, \se )$ is increasing with respect to $z \in \Se$ no greater than $\se \in \Se$ and is decreasing with respect to $z \in \Se$
no less than $\se \in \Se$.

(v) \bel & & \Pr \li \{ \wh{\bs{\se}} \geq z \mid z \ri \} \leq
\f{1}{2} + \f{C_{\mrm{BE}}}{\sq{n}} \f{ \bb{E} [ |T(X) - z |^3 ] } {
\bb{E}^{\f{3}{2}} [ |T(X) - z |^2 ] } \leq \f{1}{2} +
\f{C_{\mrm{BE}}}{\sq{n}} \f{ \bb{E}^{\f{3}{4}}
[ |T(X) - z |^4 ] } { \bb{E}^{\f{3}{2}} [ |T(X) - z |^2 ] }, \la{ineq8}\\
& & \Pr \li \{ \wh{\bs{\se}} \leq z \mid z \ri \} \leq \f{1}{2} +
\f{C_{\mrm{BE}}}{\sq{n}} \f{ \bb{E} [ |T(X) - z |^3 ] } {
\bb{E}^{\f{3}{2}} [ |T(X) - z |^2 ] } \leq \f{1}{2} +
\f{C_{\mrm{BE}}}{\sq{n}} \f{ \bb{E}^{\f{3}{4}} [ |T(X) - z |^4 ] } {
\bb{E}^{\f{3}{2}} [ |T(X) - z |^2 ] }, \la{ineq88} \eel where the
expectation is taken with $\se = z$ and $C_{\mrm{BE}}$ is the
absolute constant in the Berry-Essen inequality.

 \eeT

The proof of Theorem \ref{THm2} is given in Appendix \ref{THm2_app}.   By the assumption that $\eta (\se)$ is increasing with respect to $\se$,
it follows from statement (ii) that \[ \mscr{M} (z, \se ) = \bec \inf_{t < 0 } \bb{E} \li [ \exp \li ( n t ( \wh{\bs{\se}} -
z) \ri ) \ri ] & \tx{for} \; z \leq \se,\\
\inf_{t > 0 } \bb{E} \li [ \exp \li ( n t ( \wh{\bs{\se}} - z) \ri )
\ri ] & \tx{for} \; z \geq \se \eec
\]
This implies that the likelihood ratio bound coincides with Chernoff
bound for the exponential family.

Theorem \ref{THm2} involves the famous Berry-Essen inequality
\cite{Berry, Essen}, which asserts the following:

    Let $Y_1, Y_2, ...$ be i.i.d. samples of random variable $Y$ such that
    $\bb{E}[Y] = 0, \; \bb{E}[Y^2] > 0$, and $\bb{E}[|Y|^3] < \iy$. Also, let
$F_n$ be the cdf of ${\sum_{i = 1}^n Y_i \over {\sqrt{n \bb{E}[Y^2]}
}}$, and $\Phi$ the cdf of the standard normal distribution. Then,
there exists a positive constant $C_{\mrm{BE}}$ such that for all
$y$ and $n$,
\[
\left|F_n(y) - \Phi(y)\right| \le \f{C_{\mrm{BE}}}{\sq{n}}
{\bb{E}[|Y|^3] \over \bb{E}^{3\sh 2} [Y^2]}.
\]
A few years ago, Shevtsova \cite{Shevtsova} proved that the constant
$C_{\mrm{BE}} < 0.7056 < \f{1}{\sq{2}}$. More recently, Tyurin
\cite{Tyurin} has shown that $C_{\mrm{BE}} < 0.4785 < \f{1}{2}$.

\sect{Bounds of  Tail Probabilities}  \la{sectail}

In this section, we shall apply our general results to derive sharp
bounds for the tail probabilities of some common distributions.

\subsection{Binomial Distribution}

The probability mass function of a Bernoulli random variable, $X$,
of mean value $p \in (0, 1)$ is given by
\[
f(x, p) \equiv \Pr\{ X = x \mid p \} = p^x (1 - p)^{1 - x} = h (x)
\exp \left( \eta (p) T(x) - A(p) \right), \qqu x \in \{0, 1\}
\]
where
\[
T(x) = x,  \qqu h(x) = 1, \qqu  \eta (p) = \ln {p \over 1-p}, \qqu
A(p) = \ln \f{1}{1-p}.
\]
Since $\f{ d A ( p ) }{ d p} = \lm \f{ d \eta ( p ) }{d p}$ holds,
making use of Theorem \ref{THm2}, we have the following results.

\begin{corollary} \la{BER}
  Let $X_1, \cd, X_n$ be i.i.d. samples of Bernoulli random variable $X$ of mean value $p
\in (0, 1)$. Define $\mscr{M} ( z, p) = z \ln \f{p}{z} + (1 - z) \ln
\f{ 1 - p }{ 1 - z}$ for $z \in (0, 1)$ and $p \in (0, 1)$.  Then,
\bee &  & \Pr \li \{ \sum_{i = 1}^n X_i \geq nz \ri \} \leq \li (
\f{1}{2} + \vDe
\ri ) \exp ( n \mscr{M} ( z, p) ) \qqu \tx{for $z \in (p, 1)$},\\
&  & \Pr \li \{ \sum_{i = 1}^n X_i \leq nz \ri \} \leq \li (
\f{1}{2} + \vDe \ri ) \exp ( n \mscr{M} ( z, p) ) \qqu \tx{for $z
\in (0, p)$}, \eee where
\[
\vDe = \min \li \{ \f{1}{2}, \f{C_{\mrm{BE}} [z^2 + (1 - z)^2] }{
\sq{ n z (1 - z) } } \ri \}.
\]

\end{corollary}

An important application of Corollary \ref{BER} can be found in the
determination of sample size for estimating binomial parameters. Let
$X_1, X_2, \cd$ be i.i.d. samples of Bernoulli random $X$ such that
$\Pr \{ X = 1 \} = 1 - \Pr \{ X = 0 \} = p \in (0,1)$.  Define $\wh{
\bs{p} }_n = \f{\sum_{i = 1}^n X_i } {n}$. A classical problem in
probability and statistics theory is as follows:

Let $\vep \in (0,1)$ and $\de \in (0, 1)$ be the margin of absolute
error and the confidence parameter respectively.   How large $n$ is
sufficient to ensue \be \la{goal8} \Pr \{ | \wh{\bs{p}}_n - p | <
\vep \}
> 1- \de
\ee for any $p \in (0, 1)$?    The best explicit bound so far is the
well-known Chernoff-Hoeffding bound which asserts that (\ref{goal8})
is guaranteed for any $p \in (0, 1)$ provided that \be \la{CB}
 n > \f{1}{2 \vep^2} \ln
\f{2}{\de}. \ee

By virtue of Corollary \ref{BER}, we have obtained better explicit
sample size bound as follows.

\beT \la{ChenA}
Let $0 < \vep < \f{3}{4}$ and {\small $ 0 < \de < 2
\exp \li (
 - \f{ 9 \ln 2 } { (3 - 4 \vep)^2  } \ri )$}.
Then, $\Pr \{ | \wh{\bs{p}} - p | < \vep \mid p \}
> 1- \de$  for any $p \in (0, 1)$ provided that \be
\la{Chenformula} n > \f{1}{2 \vep^2} \ln  \f{1 + \zeta}{\de}, \ee
where {\small \[ \zeta =  \f{4 C_{\mrm{BE}} } { \sq{
 \li [ 1 - \li ( \f{4 \vep}{3} +  \sq{ \f{ \ln 2 }
{ \ln \f{2}{\de} }}  \ri )^2 \ri ] \f{\ln \f{1}{\de} } { 2 \vep^2} }
}.
 \]}
\eeT

\begin{figure}
\centering
\includegraphics[height=8cm]{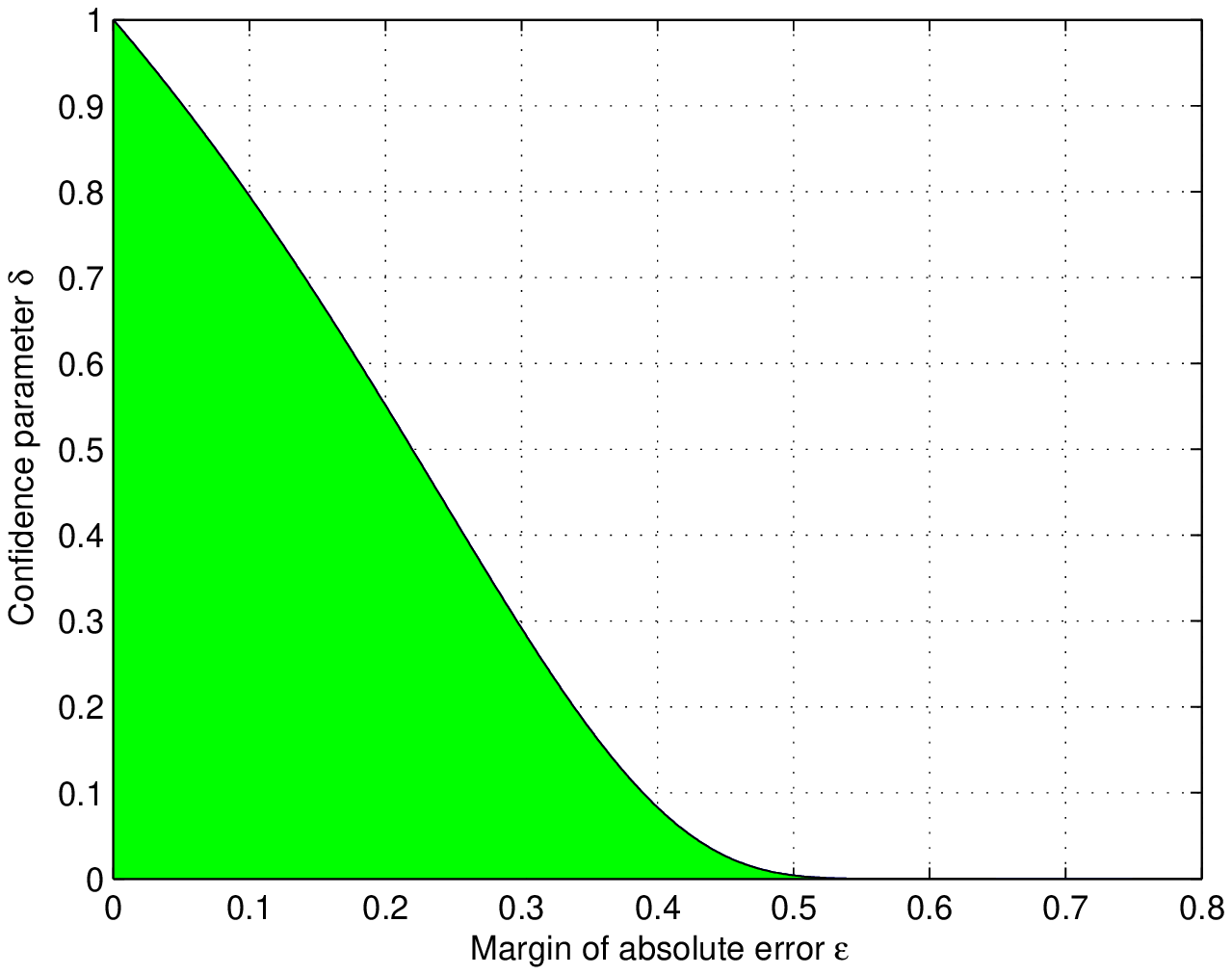}
\caption{Region of $(\vep, \de)$}
\label{Chen_dom}       
\end{figure}

The domain of $(\vep, \de)$ for which our sample size bound
(\ref{Chenformula}) can be used is shown by Figure \ref{Chen_dom}.
Clearly, a sufficient but not necessary condition to use our formula
(\ref{Chenformula}) is $0 < \vep < \f{1}{4}, \; 0 < \de < \f{1}{4}$.

\begin{figure}
\centering
\includegraphics[height=8cm]{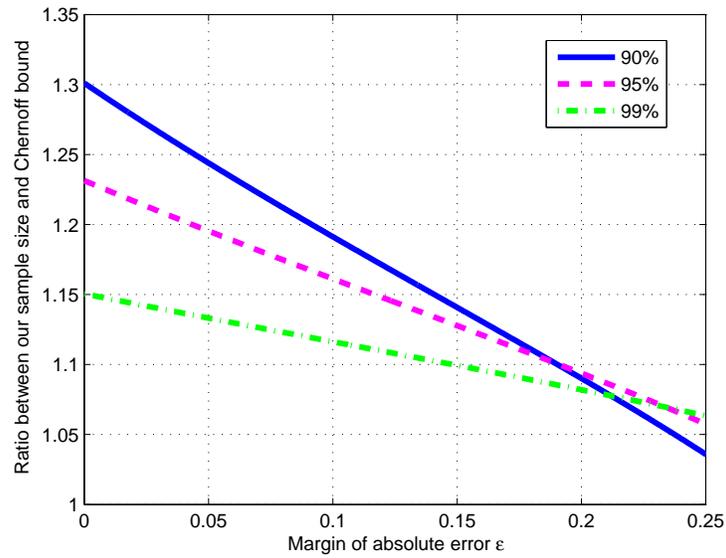}
\caption{Comparison with Chernoff bound}
\label{CHEN_CHErnoff}       
\end{figure}

The improvement of our sample size bound (\ref{Chenformula}) upon
Chernoff-Hoeffding bound (\ref{goal8}) is shown by Figure
\ref{CHEN_CHErnoff}. It can be seen that for a typical requirement
of confidence level $100 (1 - \de) \%$ (e.g., $95 \%$), the
improvement can be $20\%$ to $30\%$.

Corollary \ref{BER} is also useful for the study of inverse binomial
sampling. Let $\ga$ be a positive integer.  Define random number
$\mbf{n}$ as the minimum integer such that the summation of
$\mbf{n}$ consecutive Bernoulli random variables of common mean $p
\in (0, 1)$ is equal to $\ga$.  In other words, $\mbf{n}$ is a
random variable satisfying $\sum_{i=1}^{ \mbf{n} - 1 } X_i < \ga =
\sum_{i=1}^{ \mbf{n} } X_i$, where $X_1, X_2, \cd$ are i.i.d.
samples of Bernoulli random $X$ such that $\Pr \{ X = 1 \} = 1 - \Pr
\{ X = 0 \} = p \in (0,1)$ as mentioned earlier. This means that
$\mbf{n}$ is the least number of Bernoulli trials of success rate $p
\in (0, 1)$ to come up with $\ga$ successes.  By virtue of Corollary
\ref{BER}, we have obtained the following results.

\begin{corollary}

\la{Co2}

\bee &  & \Pr \li \{ \f{\ga}{ \mbf{n} } \leq z \ri \} \leq \li (
\f{1}{2} + \vDe \ri ) \exp \li ( \f{\ga}{z} \mscr{M} ( z, p) \ri )
\qqu \tx{for $z \in (0, p)$ such that $\f{\ga}{z}$ is an integer},\\
&  & \Pr \li \{ \f{\ga}{ \mbf{n} } \geq z \ri \} \leq \li ( \f{1}{2}
+ \vDe \ri ) \exp \li ( \f{\ga}{z} \mscr{M} ( z, p) \ri ) \qqu
\tx{for $z \in (p, 1)$ such that $\f{\ga}{z}$ is an integer}, \eee
where
\[
\vDe = \min \li \{ \f{1}{2}, \f{C_{\mrm{BE}} [z^2 + (1 - z)^2]}{
\sq{ \ga (1 - z) } }
 \ri \}.
\]
\end{corollary}

Similar to the sample size problem associated with (\ref{goal8}), it
is an important problem to estimate the binomial parameter $p$ with
a relative precision.  Specifically, consider an inverse binomial
sampling scheme as described above.  Define $\wh{\bs{p}}_\ga =
\f{\ga}{ \mbf{n}}$ as an estimator for $p$. A fundamental problem of
practical importance is stated as follows:

Let $\vep \in (0,1)$ and $\de \in (0, 1)$ be the margin of relative
error and the confidence parameter respectively.    How large $\ga$
is sufficient to ensue \be \la{revgoal8} \Pr \li \{ \li | \f{
\wh{\bs{p}}_\ga - p}{p}  \ri | < \vep \ri \}
> 1- \de
\ee for any $p \in (0, 1)$?

By virtue of Corollary \ref{Co2}, we have established the following
results regarding the above question.

\beT  The following statements (I) and (II) hold true.

(I) {\small $\Pr \li \{ \li | \f{ \wh{\bs{p}}_\ga - p}{p} \ri | <
\vep  \ri \} > 1 - \de$} for any $p \in (0, 1)$ provided that $\vep
> 0, \; 0 < \de < 1$ and \be \la{forrev0} \ga > \f{(1 + \vep)}{  (1
+ \vep) \ln (1 + \vep) - \vep } \ln \f{2}{\de}. \ee

(II) {\small $\Pr \li \{ \li | \f{ \wh{\bs{p}}_\ga  - p}{p} \ri | <
\vep  \ri \} > 1 - \de$} for any $p \in (0, 1)$ provided that $0 <
\vep < 1$, {\small \[ 0 < \de < \exp \li ( - \f{ 3 \vep^3 ( 4 +
\vep) + 4 \vep (3 + \vep) \ln 2 } { 4 (9 - 6 \vep - 2 \vep^2) }  -
\vep \sq{ \li [ \f{ 3 \vep^2 ( 4 + \vep) + 4 (3 + \vep) \ln 2 } { 4
(9 - 6 \vep - 2 \vep^2)   } \ri ]^2 + \f{  3 ( 1 + \vep)  (3 + \vep)
\ln 2 } { 2 (9 - 6 \vep - 2 \vep^2) } }  \ri )
\]}
and \be \la{forrev} \ga
> \f{(1 + \vep)}{ (1 + \vep) \ln (1 + \vep) - \vep } \ln \f{1 +
\ze}{\de}, \ee where {\small $\ze = 2 C_{\mrm{BE}} \sq{ \f{1}{m} +
\f{z}{ m - z - m z} }$} with  {\small $m = \f{ 2  }{ \vep^2 } \ln
\f{1}{\de}$} and {\small $z = 1 + \f{2 \vep}{3 + \vep} - \f{ 9 } {
( 3 + \vep )^2 } \f{\ln \f{1}{\de} } { \ln \f{2}{\de} }$}.  \eeT

\begin{figure}
\centering
\includegraphics[height=8cm]{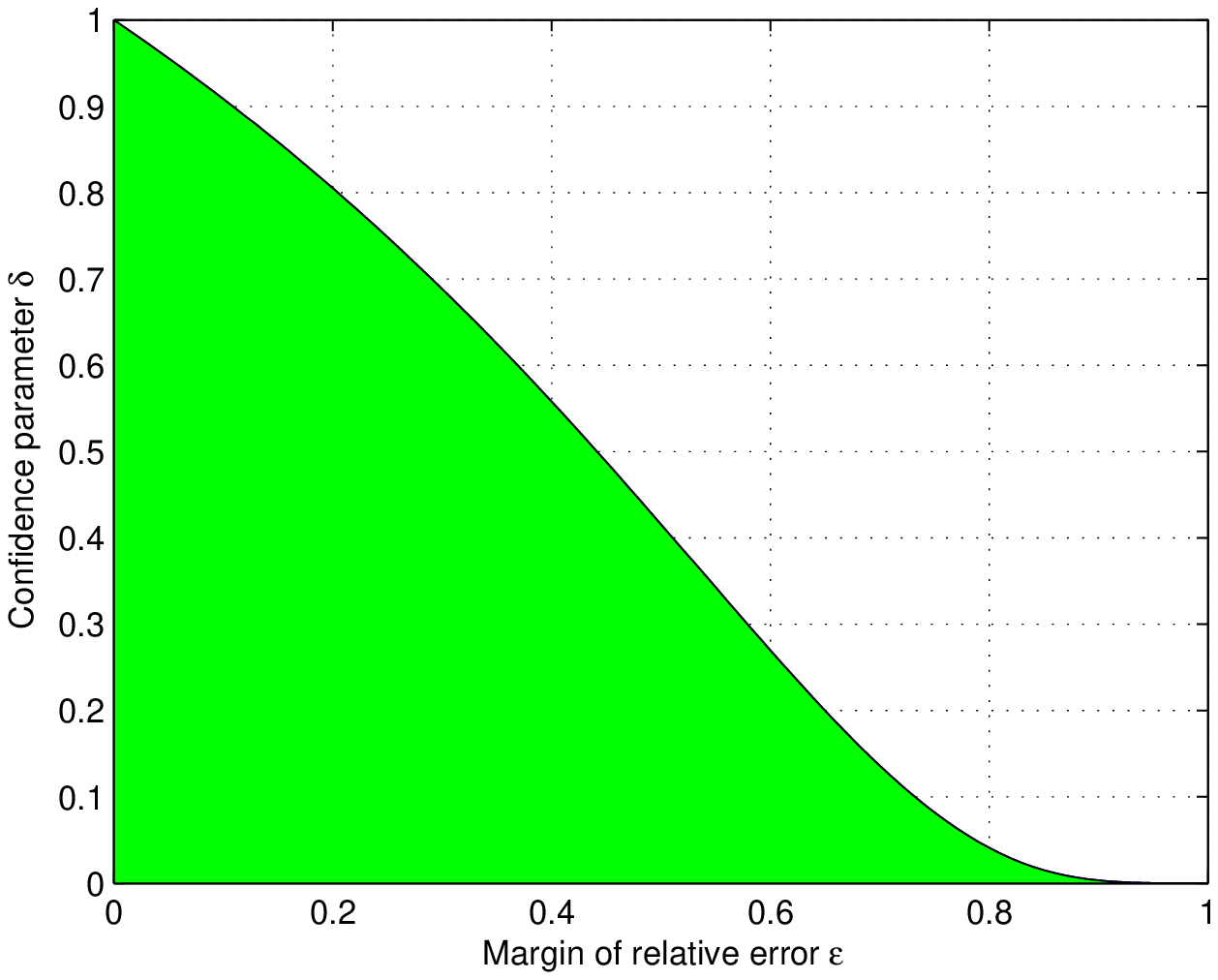}
\caption{Region of $(\vep, \de)$}
\label{Chen_dom2}       
\end{figure}

The domain of $(\vep, \de)$ for which our sample size bound
(\ref{forrev}) can be used is shown by Figure \ref{Chen_dom2}.
Clearly, a sufficient but not necessary condition to use our formula
(\ref{forrev}) is $0 < \vep < \f{3}{5}, \; 0 < \de < \f{1}{4}$.

\subsection{Negative Binomial Distribution}

The probability mass function of a negative binomial random
variable, $X$, is given by
\[
    f(x, \se) \equiv \Pr\{ X = x \mid p \} = \frac{\Gamma(x+r)}{\Ga(x + 1) \; \Gamma(r)}
    (1-p)^x p^r = h (x) \exp \left( \eta (\se) T(x) - A(\se) \right) \quad \text{for } x = 0, 1, 2,
    \dots
\]
where $r$ is a real, positive number, \[  T(x) = \f{r + x}{r}, \qqu
h(x) = \frac{\Gamma(x+r)}{\Ga(x + 1) \; \Gamma(r)}, \qqu \se =
\f{1}{p}, \qqu \eta (\se) = r \ln \li ( 1 - \f{1}{\se} \ri ), \qqu
A(\se) = r \ln (\se - 1).
\]
Since $\f{ d A ( \se ) }{ d \se } = \se \f{ d \eta ( \se ) }{d \se}$
holds, by Theorem \ref{THm2}, we have the following result.

\begin{corollary}  Let $X_1, \cd, X_n$ be i.i.d. samples of negative
binomial random variable $X$ parameterized by $\se = \f{1}{p}$.
Then,

\bee &  & \Pr \li \{ \sum_{i = 1}^n T(X_i) \geq n z \ri \} \leq \li
[  \f{p z - p}{1 - p}  \li ( \f{ z - z p}{z - 1} \ri )^z \ri ]^{n r}
\qqu
\tx{for $1 > z \geq \se = \f{1}{p}$},\\
&  & \Pr \li \{ \sum_{i = 1}^n T(X_i) \leq n z \ri \} \leq \li [
\f{p z - p}{1 - p}  \li ( \f{ z - z p}{z - 1} \ri )^z \ri ]^{n r}
\qqu \tx{for $0 < z \leq \se = \f{1}{p}$}. \eee

\end{corollary}

\subsection{Poisson Distribution}

The probability mass function of a Poisson random variable, $X$, of
mean value $\lm$ is given by
\[
f(x, \lm) \equiv \Pr \{ X = x \mid \lm \} = \f{ \lm^{x} e^{-\lm} } {
x! } = h (x) \exp \left( \eta (\lm) T(x) - A(\lm) \right), \qqu x
\in \{ 0, 1, 2, \cd \}
\]
where
\[
T(x) = x,  \qqu h(x) = \f{1}{x!}, \qqu  \eta (\lm) = \ln \lm, \qqu
A(\lm) = \lm.
\]
The moment generating function is $M(t) = \bb{E} [ e^{t X} ] = e^{-
\lm} \exp ( \lm e^t )$. Clearly, $M^\prime (t) = \lm e^t M(t)$ and
$\bb{E} [ X] = M^\prime (0) = \lm$.  It can be shown by induction
that
\[
\f{ d^{\ell + 1} M (t)  } {  d t^{\ell + 1} }  = \li ( 1 + 2^{\ell -
1} \lm e^t \ri ) \f{ d^\ell M (t) } {  d t^\ell }, \qqu \bb{E} [
X^{\ell + 1} ] = \li ( 1 + 2^{\ell - 1} \lm \ri )  \li. \f{ d^\ell M
(t) } { d t^\ell } \ri |_{t = 0} = \lm \prod_{i = 1}^\ell \li ( 1 +
2^{i - 1} \lm \ri )
\]
for $\ell = 1, 2, \cd$.  Hence, \[ \bb{E} [ |X - \lm|^2 ] = \lm,
\qqu \bb{E} [ |X - \lm|^4 ] = \sum_{i = 0}^4 \bi{4}{i}  (- \lm)^i
\bb{E} [ X^{4 - i}] = \lm ( 3 \lm^3 + 8 \lm^2 + 3 \lm + 1) \] and
\be \la{ly8} \f { \bb{E}^{\f{3}{4}} [ |X - \lm|^4 ] } {
\bb{E}^{\f{3}{2}} [ |X - \lm|^2 ] }  = \li ( 3 \lm^2 + 8 \lm + 3 +
\f{1}{\lm} \ri )^{3 \sh 4}. \ee Since $\f{ d A ( \lm ) }{ d \lm } =
\lm \f{ d \eta ( \lm ) }{d \lm}$ holds, making use of (\ref{ly8})
and Theorem 2, we have the following results.

\begin{corollary}
Let $X_1, \cd, X_n$ be i.i.d. samples of Poisson random variable $X$
of mean value $\lm$.  Then,  \bee &  & \Pr \li \{ \f{ \sum_{i = 1}^n
X_i } {n} \geq z \mid \lm \ri \} \leq \li ( \f{1}{2} + \vDe \ri )
\li ( \f{ \lm^z e^{z} } {
z^z e^{\lm} } \ri )^n \qqu \tx{for $z \geq \lm$},\\
&  & \Pr \li \{ \f{ \sum_{i = 1}^n X_i  } {n} \leq z \mid \lm  \ri
\} \leq \li ( \f{1}{2} + \vDe \ri ) \li ( \f{ \lm^z e^{z} } { z^z
e^{\lm} } \ri )^n \qqu \tx{for $0 < z \leq \lm$}, \eee where
\[
\vDe = \min \li \{ \f{1}{2}, \f{ C_{\mrm{BE}} } { \sq{n}  }  \li ( 3
z^2 + 8 z + 3 + \f{1}{z} \ri )^{ \f{3}{4}} \ri \}.
\]

\end{corollary}

\subsection{Hypergeometric Distribution}

The hypergeometric distribution can be described by the following
model.  Consider a finite population of $N$ units, of which there
are $M$ units having a certain attribute.  Draw $n$ units from the
whole population by sampling without replacement. Let $K$ denote the
number of units having the attribute found in the $n$ draws.  Then,
$K$ is a random variable possessing a hypergeometric distribution
such that
\[
\Pr \{ K = k \} = \f{ \bi{M}{k} \bi{N - M}{n - k} } { \bi{N}{n}  },
\qqu k = 0, 1, \cd, n.
\]
It can be verified that
\[
\f{ \Pr \{ K = k + 1 \mid M_1  \}  } { \Pr \{ K = k + 1 \mid M_0 \}
} \li [ \f{ \Pr \{ K = k \mid M_1  \}  } { \Pr \{ K = k \mid M_0 \}
} \ri ]^{-1}  =  \f{ (M_1 - k) ( N - M_0 - n + k + 1 ) } { (M_0 - k)
(N - M_1 - n + k + 1) } \geq 1
\]
for $M_1 \geq M_0$, which implies that the hypergeometric
distribution possesses the MLRP.  Consequently, applying Theorem
\ref{THm1}, we have the following results.

\begin{corollary}

\la{Finite}

Let $\wh{M} = \wh{M} (k)$ be a function of $k \in I_K$, which takes
values in $\{m \in \bb{Z}: k \leq m \leq N \}$. Then, \bee &  & \Pr
\{ K \leq k \mid M \} \leq \f{ \bi{ M }{k} \bi{N - M } {n - k} } {
\bi{ \wh{M} }{k} \bi{N - \wh{M} } {n - k}  }
\qu \tx{for $k \in I_K$ such that $\wh{M} (k) \leq M$},\\
&   & \Pr \{ K \geq k \mid M  \} \leq \f{ \bi{ M }{k} \bi{N - M } {n
- k} } { \bi{ \wh{M} }{k} \bi{N - \wh{M} } {n - k}  } \qu \tx{for $k
\in I_K$ such that $\wh{M} (k) \geq M$}. \eee

\end{corollary}

Actually, a specialized version of the inequalities in Corollary \ref{Finite} had been used in the $15$-th version of our paper \cite{Chen1}
published in arXiv on August 6, 2010 for developing multistage sampling schemes for estimating population proportion $p$. Moreover, the
specialized inequalities had been used in the $20$-th version of our paper \cite{Chen12} published in arXiv on August 7, 2010  for developing
multistage testing plans for hypotheses regarding $p$.

\subsection{Hypergeometric Waiting-Time Distribution}

The hypergeometric waiting-time distribution can be described by the
following model.  Consider a finite population of $N$ units, of
which there are $M$ units having a certain attribute.   Continue
sampling until $r$ units of certain attribute is observed or the
whole population is checked.  Let $\bs{n}$ be the number of units
checked when the sampling is stopped.  Clearly, in the case of $r >
M$, it must be true that $\Pr \{ \bs{n} = N \} = 1$, since the whole
population is checked. In the case of $r \leq M$, the random
variable $\bs{n}$ has a hypergeometric waiting-time distribution
such that \bee
 \Pr \{ \bs{n} = n \mid M \}  =  \f{ \bi{n - 1}{r - 1} \bi{N - n}{M - r} } {
\bi{N}{M} }  & \tx{  } \eee
 for $r \leq M$ and $r \leq n \leq N$.  It can be shown that
 \bee \f{ \Pr \{ \bs{n} = n + 1 \mid M_1 \} } {  \Pr \{
\bs{n} = n + 1 \mid M_0 \} } \li [ \f{ \Pr \{ \bs{n} = n  \mid M_1
\} } { \Pr \{ \bs{n} = n \mid M_0 \} } \ri ]^{-1} = \f{ N - n - M_1
+ r }{ N - n - M_0 + r } \geq 1 \eee for $M_0 \leq M_1$, which
implies that the hypergeometric waiting-time distribution possesses
the MLRP.  Hence, by virtue of Theorem \ref{THm1}, we have the
following results.

\begin{corollary}

Let $\wh{M} = \wh{M} (n)$ be a function of $n \in I_{\bs{n}}$, which
takes values in $\{ m \in \bb{Z}: r \leq m \leq N \}$.  Then, \bee &
& \Pr \{ \bs{n} \leq n \mid M \} \leq \f{ \bi{N}{\wh{M}} \bi{N -
n}{M - r} } { \bi{N}{M} \bi{N - n}{\wh{M} - r} } = \f{ \bi{ M }{r}
\bi{N - M } {n - r} } { \bi{ \wh{M} }{r} \bi{N - \wh{M} } {n - r}  }
\qu \tx{for $n \in I_{\bs{n}}$ such that $\wh{M} (n) \geq M$},\\
&   & \Pr \{ \bs{n} \geq n \mid M \} \leq \f{  \bi{N}{\wh{M}} \bi{N
- n}{M - r} } { \bi{N}{M} \bi{N - n}{\wh{M} - r} } = \f{ \bi{ M }{r}
\bi{N - M } {n - r} } { \bi{ \wh{M} }{r} \bi{N - \wh{M} } {n - r}  }
\qu \tx{for $n \in I_{\bs{n}}$ such that $\wh{M} (n) \leq M$}. \eee

\end{corollary}

\subsection{Normal Distribution}

The probability density function of a Gaussian random variable, $X$,
with mean $\mu$ and variance $\si^2$ is given by
\[
    f (x;\mu) = \frac{1}{\sqrt{2 \pi} \sigma} \exp \li (
    -\f{|x-\mu|^2}{2\sigma^2} \ri )
    = h (x) \exp \left( \eta (\se) T(x) - A(\se) \right).
\]
where
\[
T(x) = \f{x}{\sigma}, \qu  h(x) = \f{1}{\sqrt{2\pi} \sigma} \exp \li
( - \f{x^2}{2\sigma^2} \ri ), \qu \se = \f{\mu}{\si}, \qu A(\se) =
\f{\se^2}{2}, \qu \eta(\se) = \se.
\]
Since $\f{ d A ( \se ) }{ d \se } = \se \f{ d \eta ( \se ) }{d \se}$
holds, by Theorem 2, we have the following results.
\begin{corollary}
\la{Gauss}  Let $X_1, \cd, X_n$ be i.i.d. samples of Gaussian random
variable $X$ of mean $\mu$ and variance $\si^2$.  Then,  \bee & &
\Pr \li \{ \f{ \sum_{i = 1}^n X_i}{n} \leq z \ri \} < \f{1}{2} \exp
\li ( - \f{(z - \mu)^2}{2 \si^2} \ri ) \qqu
\tx{for $z \leq \mu$},\\
&  & \Pr \li \{  \f{ \sum_{i = 1}^n X_i}{n} \geq z \ri \} < \f{1}{2}
\exp \li ( - \f{(z - \mu)^2}{2 \si^2} \ri ) \qqu \tx{for $z \geq
\mu$}. \eee
\end{corollary}

It should be noted that the inequalities in Corollary \ref{Gauss}
may be shown by using other methods.  However, the factor $\f{1}{2}$
cannot be obtained by using Chernoff bounds.

\subsection{Gamma Distribution}

In probability theory and statistics, a random variable $X$ is said
to have a gamma distribution if its density function is of the form
\[
f(x) = \frac{x^{k - 1}} { \Gamma(k) \se ^{ k }    }  \exp \li ( -
\frac{x}{\se} \ri ) = h (x) \exp \left( \eta (\se) T(x) - A(\se)
\right) \;\;\; {\rm for} \;\;\; 0 < x < \infty
\]
where $\se > 0, \; k > 0$ are referred to as the scale parameter and
shape parameter respectively, and
\[ h(x) = \f{ x^{k - 1} } { \Ga(k)  }, \qqu T (x) = \f{x}{k}, \qqu
\eta (\se) = - \f{k}{\se}, \qqu A ( \se ) = k \ln \se.
\]
The moment generating function of $X$ is $M(t) = \bb{E} [ e^{t X} ]
= (1 - \se t)^{-k}$ for $t < \f{1}{\se}$.  It can be shown by
induction that
\[
\f{ d^{\ell + 1} M (t)  } {  d t^{\ell + 1} }  = \f{( k + \ell)
\se}{1 - \se t} \f{ d^\ell M (t) } {  d t^\ell }, \qqu \bb{E} [
X^{\ell + 1} ] = (k + \ell) \se \li. \f{ d^\ell M (t) } {  d t^\ell
} \ri |_{t = 0} = \se^{\ell + 1} \prod_{i = 0}^\ell (k + i)
\]
for $\ell = 0, 1, 2, \cd$. Therefore, \[ \bb{E} [ | X - k \se |^2 ]
= k \se^2, \qqu  \bb{E} [ | X - k \se |^4 ] = \sum_{i = 0}^4
\bi{4}{i}  (- k \se)^i \bb{E} [ X^{4 - i}] = 3 k (k + 2) \se^4, \]
and \be \la{lya} \f{  \bb{E}^{ \f{3}{4} } [ | X - k \se |^4 ]  } {
\bb{E}^{ \f{3}{2} } [ | X - k \se |^2 ]  } = \li ( 3 + \f{6}{k} \ri
)^{ \f{3}{4} }. \ee Since $\f{ d A ( \se ) }{ d \se } = \se \f{ d
\eta ( \se ) }{d \se}$ holds, making use of (\ref{lya}) and Theorem
2, we have

\begin{corollary}
\la{Gammac}  Let $X_1, \cd, X_n$ be i.i.d. samples of Gamma random
variable $X$ of shape parameter $k$ and scale parameter $\se$. Then,
\bee & & \Pr \li \{ \frac{1}{n} \sum_{i=1}^n X_i \geq \ro k \se \ri
\} \leq \li ( \f{1}{2} + \vDe \ri ) \li [ \ro  \exp \li ( 1 - \ro
\ri ) \ri ]^{k n} \qqu \tx{for $\ro \geq 1$},\\
&  & \Pr \li \{ \frac{1}{n} \sum_{i=1}^n X_i \leq \ro k \se \ri \}
\leq \li ( \f{1}{2} + \vDe \ri ) \li [ \ro  \exp \li ( 1 - \ro \ri )
\ri ]^{k n} \qqu \tx{for $0 < \ro \leq 1$}, \eee where
\[
\vDe = \min \li \{ \f{1}{2},  \li ( 3 + \f{6}{k} \ri )^{ \f{3}{4} }
\f{ C_{\mrm{BE}} } { \sq{n} }  \ri \}.
\]

\end{corollary}

It should noted that the chi-square distribution of $k$ degrees of
freedom is a special case of the Gamma distribution with shape
parameter $\f{k}{2}$ and scale parameter $2$.  The exponential
distribution of mean $\se$ is also a special case of the Gamma
distribution with shape parameter $1$ and scale parameter $\se$. If
the shape parameter $k$ is an integer, then the Gamma distribution
represents an Erlang distribution.  Therefore, the bounds in
Corollary \ref{Gammac} can be used for those distributions.

Let $\wh{\bs{\se}} = \f{ \sum_{i = 1}^n X_i  } { k n }$.  In order
to find the sample size $n$ such that $\Pr \li \{  \li |
\wh{\bs{\se}} - \se \ri | < \vep  \se \ri \} > 1 - \de$, we have
established the following result.

\beT Let $\vep > 0$ and $0 < \de < 1$.  Then,  $\Pr \li \{  \li |
\wh{\bs{\se}} - \se \ri | < \vep \se \ri \} > 1 - \de$ if $n >
\f{\ln \f{1 + \ze}{\de}} { k [\vep + \ln(1 + \vep) ] }$, where
\[
\ze = 2 C_{\mrm{BE}} \li ( 3 + \f{6}{k} \ri )^{ \f{3}{4} } \sq{\f{ k
[ \vep + \ln(1 + \vep) ] }{\ln \f{1}{\de}} }.
\]
\eeT

\subsection{Student's $t$-Distribution}

If the random variable $X$ has a density function of the form
\[
f(x) = \f{ \Ga ( \f{n + 1}{2}  ) } { \sq{n \pi} \Ga (\f{n}{2} )  ( 1
+ \f{x^2}{n}  )^{(n + 1) \sh 2} }, \qqu \tx{for} \qu - \iy < x <
\iy,
\]
then the variable $X$ is said to posses a Student's $t$-distribution
with $n$ degrees of freedom.

Now, we want to bound the tail probabilities of the distribution of
$X$. Define $Y = \se |X|$,  where $\se$ is a positive number.  Then,
$Y$ is a random variable parameterized by $\se$.  For any real
number $t$,
\[
\Pr \{ Y \leq t \} =  \Pr \li \{  |X| \leq \f{t}{\se} \ri \} = 2
\int_0^{ \f{t}{\se} } f(x) dx - 1.
\]
By differentiation, we obtain the probability density function of
$Y$ as $f_{Y} (t, \se) = \f{2}{\se}  f \li ( \f{t}{\se} \ri )$. Note
that, for $\se_0 < \se_1$,  \bee \f{ f_{Y} (t, \se_1) } { f_{Y} (t,
\se_0) }  = \f{\se_0}{\se_1} \f{ f \li ( \f{t}{\se_1} \ri ) } { f
\li ( \f{t}{\se_0} \ri ) } = \f{\se_0}{\se_1} \li ( 1 + \f{ \f{1}{
\se_0^2} - \f{1}{ \se_1^2}  } { \f{n}{t^2} + \f{1}{\se_1^2} } \ri
)^{ (n + 1) \sh 2 }, \eee which is monotonically increasing with
respect to $t \in I_{Y}$. This implies that the likelihood ratio
$\f{ f_{Y} (t, \se_1) } { f_{Y} (t, \se_0) }$ is monotonically
increasing with respect to $Y$. Therefore, by Theorem \ref{THm1},
\bee \Pr \{ |X| \geq x \}  = \Pr \{ Y \geq x \se \} \leq \f{
\f{2}{\se} f \li ( \f{x \se}{\se} \ri )  } { \f{2}{x \se} f \li (
\f{x \se}{x \se} \ri ) } = \f{ x f(x) } { f(1) } = x \li ( \f{ n +
1} {n + x^2} \ri )^{ (n+1) \sh 2} \eee for $x \geq 1$.  Similarly,
\bee \Pr \{ |X| \leq x \}  = \Pr \{ Y \leq x \se \} \leq \f{
\f{2}{\se}  f \li ( \f{x \se}{\se} \ri )  } { \f{2}{x \se} f \li (
\f{x \se}{x \se} \ri ) } = \f{ x f(x) } { f(1) } = x \li ( \f{ n +
1} {n + x^2} \ri )^{ (n+1) \sh 2} \eee for $0 < x \leq 1$.  By
differentiation, we can show that the upper bound of the tail
probabilities is unimodal with respect to $x$. In summary, we have
the following results.

\begin{corollary}

Suppose $X$ possesses a Student's $t$-distribution with $n$ degrees
of freedom.  Then, \bee  &   &  \Pr \{ |X| \geq x \} \leq  x \li (
\f{ n + 1} {n + x^2}
\ri )^{ (n+1) \sh 2} \qqu \tx{for $x \geq 1$},\\
&   &  \Pr \{ |X| \leq x \} \leq  x \li ( \f{ n + 1} {n + x^2} \ri
)^{ (n+1) \sh 2} \qqu \tx{for $0 \leq x \leq 1$}, \eee where the
upper bound of the tail probabilities is monotonically increasing
with respect to $x \in (0, 1)$ and monotonically decreasing with
respect to $x \in (1, \iy)$.

\end{corollary}

\subsection{Snedecor's $F$-Distribution}

If the random variable $X$ has a density function of the form
\[
f(x) = \f{ \Ga ( \f{n + m}{2}  )  ( \f{m}{n} )^{m \sh 2}  x ^{ (m -
2) \sh 2} } { \Ga ( \f{m}{2} ) \Ga (\f{n}{2} )  ( 1 +  \f{m}{n} x
)^{(n + m) \sh 2} }, \qqu \tx{for} \qu 0 < x < \iy,
\]
then the variable $X$ is said to posses an $F$-distribution with $m$
and $n$ degrees of freedom.

Now, we want to bound the tail probabilities of the distribution of
$X$. Define $Y = \se X$,  where $\se$ is a positive number.  Then,
$Y$ is a random variable parameterized by $\se$.  For any real
number $t$,
\[
\Pr \{ Y \leq t \} =  \Pr \li \{  X \leq \f{t}{\se} \ri \} =
\int_0^{ \f{t}{\se} } f(x) dx.
\]
By differentiation, we obtain the probability density function of
$Y$ as $f_{Y} (t, \se) = \f{1}{\se}  f \li ( \f{t}{\se} \ri )$. Note
that, for $\se_0 < \se_1$,  \bee \f{ f_{Y} (t, \se_1) } { f_{Y} (t,
\se_0) }  = \f{\se_0}{\se_1} \f{ f \li ( \f{t}{\se_1} \ri ) } { f
\li ( \f{t}{\se_0} \ri ) } = \li ( \f{\se_0}{\se_1} \ri )^{m \sh 2}
\li ( 1 + \f{ \f{1}{\se_0} - \f{1}{\se_1}  } { \f{n}{mt} +
\f{1}{\se_1} } \ri )^{ (n + m) \sh 2 }, \eee which is monotonically
increasing with respect to $t \in I_{Y}$. This implies that the
likelihood ratio $\f{ f_{Y} (t, \se_1) } { f_{Y} (t, \se_0) }$ is
monotonically increasing with respect to $Y$. Therefore, by Theorem
1, \bee \Pr \{ X \geq x \}  = \Pr \{ Y \geq x \se \} \leq \f{
\f{1}{\se}  f \li ( \f{x \se}{\se} \ri )  } { \f{1}{x \se} f \li (
\f{x \se}{x \se} \ri ) } = \f{ x f(x) } { f(1) } = x^{m \sh 2} \li (
\f{ n + m} {n + m x } \ri )^{ (m + n) \sh 2} \eee for $x \geq 1$.
Similarly, \bee \Pr \{ X \leq x \}  =  \Pr \{ Y \leq x \se \} \leq
\f{ \f{1}{\se}  f \li ( \f{x \se}{\se} \ri )  } { \f{1}{x \se} f \li
( \f{x \se}{x \se} \ri ) } = \f{ x f(x) } { f(1) } = x^{m \sh 2} \li
( \f{ n + m} {n + m x } \ri )^{ (m + n) \sh 2} \eee for $0 < x \leq
1$.  By differentiation, we can show that the upper bound of the
tail probabilities is unimodal with respect to $x$.
 Formally, we state the results as follows.

\begin{corollary}

Suppose $X$ possesses an $F$-distribution with $m$ and $n$ degrees
of freedom.  Then,
\bee &  & \Pr \{ X \geq x \}  \leq x^{m \sh 2}
\li ( \f{ n + m} {n +
m x } \ri )^{ (m + n) \sh 2} \qqu \tx{for $x \geq 1$}, \\
&  & \Pr \{ X \leq x \}  \leq  x^{m \sh 2} \li ( \f{ n + m} {n + m
x} \ri )^{ (m + n) \sh 2} \qqu \tx{for $0 < x \leq 1$}, \eee  where
the upper bound of the tail probabilities is monotonically
increasing with respect to $x \in (0, 1)$ and monotonically
decreasing with respect to $x \in (1, \iy)$.

\end{corollary}

\section{Using Probabilistic Inequalities for Parameter Estimation}  \la{secest}

In this section, we shall explore the general applications of the
probabilistic inequalities for parameter estimation.

\subsection{Interval Estimation}

From Theorem \ref{THm1}, it can be seen that,  for a large class of
distributions, the likelihood ratio bounds of the cumulative
distribution function and complementary cumulative distribution of
random variable $\bs{\varphi}$ are partially monotone. Such
monotonicity can be explored for the interval estimation of the
underlying parameter $\se$. In this direction, we have developed a
method for constructing a confidence interval for $\se$ as follows.

\beT \la{LRB_CI}

Let $\bs{\varphi}$ be a random variable possessing a distribution
determined by parameter $\se \in \Se$. Let $I_{\bs{\varphi}}$ denote
the support of $\bs{\varphi}$. Let $\bb{F} (., .)$ and $\bb{G}
(.,.)$ be bivariate functions possessing the following properties:

(i) $\bb{F} (z, \vse)$ is non-increasing with respect to $\vse$ no
less than $z \in I_{\bs{\varphi}}$;

(ii) $\bb{G} (z, \vse)$ is non-decreasing with respect to $\vse$ no
greater than $z \in I_{\bs{\varphi}}$;

(iii) \bee &  & \Pr \{ \bs{\varphi} \leq z \mid \se \} \leq \bb{F}
(z, \se) \qu \tx{for $z$ no greater than } \; \se \in
\Se,\\
&  & \Pr \{ \bs{\varphi} \geq z \mid \se \} \leq \bb{G} (z, \se) \qu
\tx{for $z$ no less than} \; \se \in \Se. \eee Let $\de \in (0, 1)$.
Define confidence limits $L (\bs{\varphi}, \de)$ and $U
(\bs{\varphi}, \de)$ as functions of $\bs{\varphi}$ and $\de$ such
that $\{ \bb{F} ( \bs{\varphi}, U (\bs{\varphi}, \de) ) \leq
\f{\de}{2}, \; \bb{G} ( \bs{\varphi}, L (\bs{\varphi}, \de) ) \leq
\f{\de}{2},  \; L (\bs{\varphi}, \de) \leq \bs{\varphi} \leq U
(\bs{\varphi}, \de) \}$ is a sure event. Then, $\Pr \{ L
(\bs{\varphi}, \de) \leq \se \leq U (\bs{\varphi}, \de) \mid \se \}
\geq 1 - \de$ for any $\se \in \Se$. \eeT

See Appendix \ref{LRB_CI_app} for a proof.  By the monotonicity of
$\bb{F}(z, \se)$ and $\bb{G}(z, \se)$ with respect to $\se$, we can
obtain the lower and upper confidence limits $L(\bs{\varphi}, \de)$
and $U(\bs{\varphi}, \de)$ by a bisection approach.  In the context
of Theorem \ref{THm1}, $\bb{F}(z, \se)$ and $\bb{G}(z, \se)$ have
the same expression $\mscr{M} (z, \se)$.

\subsection{Asymptotically Tight Bound of  Sample Size}

Clearly, the likelihood ratio bound may be applied to the
determination of sample size for parameter estimation.  Since the
likelihood ratio bound coincides with Chernoff bound for the
exponential family, it is interesting to investigate the sample size
issue in connection with Chernoff bound.

Let a population be denoted by a random variable $X$. Let $\mu$ be
the mean of $X$. Suppose that the distribution of $X$ is
parameterized by $\mu$.  Suppose that the moment generating function
 $\bb{E} [ e^{t X}]$ exists for any real number $t$.  Let $\ovl{X}_n = \f{ \sum_{i = 1}^n
X_i } { n }$,  where $X_1, \cd, X_n$ are i.i.d. samples of random
variable $X$. Chernoff bound asserts that
\[
\Pr \{ \ovl{X}_n \leq \mu - \vep \} \leq  [ \mcal{F} (\mu - \vep,
\mu)  ]^n,
\]
\[
\Pr \{ \ovl{X}_n \geq \mu + \vep \} \leq  [ \mcal{G} (\mu + \vep,
\mu)  ]^n
\]
where \[ \mcal{F} (\mu - \vep, \mu) = \inf_{t < 0} \bb{E} [ e^{t (X
- \mu + \vep) } ], \qqu \mcal{G} (\mu + \vep, \mu) = \inf_{t > 0}
\bb{E} [ e^{t (X - \mu - \vep ) } ].
\]
Let $\vep > 0$ be a pre-specified margin of absolute error.  Let
$\de > 0$ be a pre-specified confidence parameter. It is a
ubiquitous problem to estimate $\mu$ by its empirical mean
$\ovl{X}_n$ such that
\[
\Pr \{ | \ovl{X}_n - \mu | < \vep  \} > 1 - \de.
\]
To guarantee the above requirement, it suffices to choose the sample
size $n$ greater than
\[
N_c (\de) \DEF \max \li \{  \f{ \ln \f{\de}{2} } {  \ln \mcal{F}
(\mu - \vep, \mu) }, \;  \f{ \ln \f{\de}{2} } {  \ln \mcal{G} (\mu +
\vep, \mu) } \ri \}.
\]
It is of theoretical and practical importance to know tightness of
such sample size bound. Let $N_{\mrm{a}} (\de)$ be the minimum
sample size $n$ to guarantee $\Pr \{ | \ovl{X} - \mu | < \vep \} > 1
- \de$.  We discover the following interesting result.

\beT  \la{AST_Bound}
\[
\lim_{\de \to 0} \f{ N_{\mrm{c}}(\de) } { N_{\mrm{a}} (\de) } = 1.
\]
\eeT

See Appendix \ref{AST_Bound_App} for a proof.  This theorem implies
that, for high confidence estimation (i.e., small $\de$), the sample
size bound $N_c (\de)$ can be quite tight.

\section{Conclusion} \la{secon}

In this paper, we have opened a new avenue for deriving
probabilistic inequalities. Especially, we have established a
fundamental connection between monotone likelihood ratio and tail
probabilities. A unified theory has been developed  for bounding the
tail probabilities of the exponential family of distributions.
Simple and sharp bounds are obtained for some other important
distributions.

\appendix

\sect{Proof of Theorem \ref{THm1} } \la{THm1_app}

To prove inequalities (\ref{bound1}) and (\ref{bound2}), we shall focus on the case that $X_1, \cd, X_n$ are discrete random variables. First,
we need to establish (\ref{bound1}).  For $z \in \mscr{Z}$ such that $\vse(z)$ is no less than $\se$, the inequality (\ref{bound1}) is trivially
true if $\mscr{M} (z, \se)$ is not bounded.   It remains to consider the case that $\mscr{M} (z, \se)$ is bounded.  By the MLRP assumption, for
$z \in \mscr{Z}$ such that $\vse(z)$ is no less than $\se$, the likelihood ratio $\Lambda(y, \se, \vse(z))$ is non-decreasing with respect to $y
\in \mscr{Z}$. In other words, the likelihood ratio $\Lambda(y, \vse(z), \se)$ is non-increasing with respect to $y \in \mscr{Z}$ provided that
$\vse(z) \geq \se$. Hence, for $z \in \mscr{Z}$ such that $\vse(z) \geq \se$, it must be true that $\Lambda ( \varphi(x_1, \cd, x_n), \vse(z),
\se ) \leq \Lambda ( z, \vse(z), \se )$ for all observation $(x_1, \cd, x_n)$ of random tuple $(X_1, \cd, X_n)$ such that $\varphi(x_1, \cd,
x_n) \geq z$. Moreover,  since $\mscr{M} (z, \se)$ is bounded, it must be true that $f_n ( x_1, \cd, x_n; \vse(z)) > 0$ for all observation
$(x_1, \cd, x_n)$ of random tuple $(X_1, \cd, X_n)$ such that $\varphi(x_1, \cd, x_n) \geq z$ and $f_n ( x_1, \cd, x_n; \se) > 0$.  It follows
that {\small \bee \Pr \{ \bs{\varphi} \geq z \mid \se \} & = & \sum_{ \varphi(x_1, \cd, x_n)
\geq z \atop{ f_n ( x_1, \cd, x_n; \se) > 0 } } f_n ( x_1, \cd, x_n; \se)\\
& = & \sum_{ \varphi(x_1, \cd, x_n) \geq z \atop{ f_n ( x_1, \cd, x_n; \se) > 0} } \f{ f_n ( x_1, \cd, x_n; \se)} { f_n ( x_1, \cd, x_n;
\vse(z)) } \times
f_n ( x_1, \cd, x_n; \vse(z))\\
& = & \sum_{ \varphi(x_1, \cd, x_n) \geq z \atop{ f_n ( x_1, \cd, x_n; \se) > 0} }  \Lambda ( \varphi(x_1, \cd, x_n), \vse(z), \se )  \times
f_n ( x_1, \cd, x_n; \vse(z))\\
& \leq & \sum_{ \varphi(x_1, \cd, x_n) \geq z }  \Lambda (
\varphi(x_1, \cd, x_n), \vse(z), \se )  \times
f_n ( x_1, \cd, x_n; \vse(z))\\
& \leq & \sum_{ \varphi(x_1, \cd, x_n) \geq z }  \Lambda ( z,
\vse(z), \se )  \times f_n ( x_1, \cd, x_n; \vse(z))\\
& = & \Lambda ( z, \vse(z), \se ) \sum_{ \varphi(x_1, \cd, x_n) \geq
z } f_n ( x_1, \cd, x_n; \vse(z))\\
&  = &  \mscr{M} (z, \se) \times \Pr \{ \bs{\varphi} \geq z \mid \vse(z) \} \leq \mscr{M} (z, \se) \eee} for $z \in \mscr{Z}$ such that
$\vse(z)$ is no less than $\se$. This establishes (\ref{bound1}).

In order to show (\ref{bound2}),  it suffices to consider the case that $\mscr{M} (z, \se)$ is bounded,  since the inequality (\ref{bound2}) is
trivially true if $\mscr{M} (z, \se)$ is not bounded.  By the MLRP assumption, for $z \in \mscr{Z}$  such that $\vse(z)$ is no greater than
$\se$, the likelihood ratio $\Lambda(y, \vse(z), \se)$ is non-decreasing with respect to $y \in \mscr{Z}$. Hence, for $z \in \mscr{Z}$ such that
$\vse(z) \leq \se$, it must be true that $\Lambda ( \varphi(x_1, \cd, x_n), \vse(z), \se ) \leq \Lambda ( z, \vse(z), \se )$ for all observation
$(x_1, \cd, x_n)$ of random tuple $(X_1, \cd, X_n)$ such that $\varphi(x_1, \cd, x_n) \leq z$.  Moreover, since $\mscr{M} (z, \se)$ is bounded,
it must be true that $f_n ( x_1, \cd, x_n; \vse(z)) > 0$ for all observation $(x_1, \cd, x_n)$ of random tuple $(X_1, \cd, X_n)$ such that
$\varphi(x_1, \cd, x_n) \leq z$ and $f_n ( x_1, \cd, x_n; \se) > 0$.  It follows that {\small \bee \Pr \{ \bs{\varphi} \leq z \mid \se \} & = &
\sum_{ \varphi(x_1, \cd, x_n)
\leq z \atop{ f_n ( x_1, \cd, x_n; \se) > 0} } f_n ( x_1, \cd, x_n; \se)\\
& = & \sum_{ \varphi(x_1, \cd, x_n) \leq z \atop{ f_n ( x_1, \cd, x_n; \se) > 0} } \f{ f_n ( x_1, \cd, x_n; \se)} { f_n ( x_1, \cd, x_n;
\vse(z)) } \times
f_n ( x_1, \cd, x_n; \vse(z))\\
& = & \sum_{ \varphi(x_1, \cd, x_n) \leq z \atop{ f_n ( x_1, \cd, x_n; \se) > 0} }  \Lambda ( \varphi(x_1, \cd, x_n), \vse(z), \se )  \times
f_n ( x_1, \cd, x_n; \vse(z))\\
& \leq & \sum_{ \varphi(x_1, \cd, x_n) \leq z }  \Lambda (
\varphi(x_1, \cd, x_n), \vse(z), \se )  \times
f_n ( x_1, \cd, x_n; \vse(z))\\
& \leq & \sum_{ \varphi(x_1, \cd, x_n) \leq z }  \Lambda ( z,
\vse(z), \se )  \times f_n ( x_1, \cd, x_n; \vse(z))\\
& = & \Lambda ( z, \vse(z), \se ) \sum_{ \varphi(x_1, \cd, x_n) \leq
z } f_n ( x_1, \cd, x_n; \vse(z))\\
&  = &  \mscr{M} (z, \se) \times \Pr \{ \bs{\varphi} \leq z \mid \vse(z) \} \leq \mscr{M} (z, \se) \eee} for $z \in \mscr{Z}$ such that
$\vse(z)$ is no greater than $\se$.  This proves (\ref{bound2}).

The proof of inequalities (\ref{bound1}) and (\ref{bound2}) for the
case that $X_1, \cd, X_n$ are continuous variables can be completed
by replacing the summation of probability mass functions with
integration of probability density functions.  It remains to show
statements (i), (ii) and (iii).

Clearly, statement (i) is a direct consequence of assumptions (a),
(b) and the definition of $\varphi(.)$.  The monotonicity of
$\mscr{M} (z, \se)$ with respect to $\se$ as described by statement
(ii) of the theorem can be established as follows.  To show
$\mscr{M} (z, \se_2) \leq \mscr{M} (z, \se_1)$ for $\se_2
> \se_1 \geq z$, note that
\[
\mscr{M} (z, \se_2) = \f{ g (z, \se_2) } { g (z, z) } \leq \f{ g (z,
\se_1) } { g (z, z) } = \mscr{M} (z, \se_1),
\]
where the inequality is due to the assumption that $g(z, \se)$ is
non-increasing with respect to $\se$ no less than $z$.  On the other
hand,  to show $\mscr{M} (z, \se_1) \leq \mscr{M} (z, \se_2)$ for
$\se_1 < \se_2 \leq z$,  note that
\[
\mscr{M} (z, \se_1) = \f{ g (z, \se_1) } { g (z, z) } \leq \f{ g (z,
\se_2) } { g (z, z) } = \mscr{M} (z, \se_2), \]  where the
inequality is due to  the assumption that $g(z, \se)$ is
non-decreasing with respect to $\se$ no greater than $z$. This
justifies statement (ii) of the theorem.

Finally, consider the monotonicity of $\mscr{M} (z, \se)$ with
respect to $z$ as described by statement (iii) of the theorem. To
show $\mscr{M} (z_2, \se) \leq \mscr{M} (z_1, \se)$ for $z_2 > z_1
\geq \se$,  notice that
\[
\mscr{M} (z_2, \se) = \f{ g (z_2, \se) } { g (z_2, z_2) } \leq \f{ g
(z_2, \se) } { g (z_2, z_1) } \leq \f{ g (z_1, \se) } { g (z_1, z_1)
} = \mscr{M} (z_1, \se),
\]
where the first inequality is due to the assumption that $g(z, \se)$
is non-decreasing with respect to $\se$ no greater than $z$ and the
second one is due to the assumption that $\Lambda(z, \se_0, \se_1) =
\f{ g(z, \se_1) }{ g(z, \se_0) }$ is non-decreasing with respect to
$z$ provided that $\se_0 \leq \se_1$. On the other side, to show
$\mscr{M} (z_2, \se) \geq \mscr{M} (z_1, \se)$ for $z_1 < z_2 \leq
\se$, it suffices to observe that
\[
\mscr{M} (z_1, \se) = \f{ g (z_1, \se) } { g (z_1, z_1) } \leq \f{ g
(z_1, \se) } { g (z_1, z_2) } \leq \f{ g (z_2, \se) } { g (z_2, z_2)
} = \mscr{M} (z_2, \se),
\]
where the first inequality is due to the assumption that $g(z, \se)$
is non-increasing with respect to $\se$ no less than $z$ and the
second one is due to the assumption that $\Lambda(z, \se_0, \se_1) =
\f{ g(z, \se_1) }{ g(z, \se_0) }$ is non-decreasing with respect to
$z$ provided that $\se_0 \leq \se_1$.  Statement (iii) of the
theorem is thus proved.

\sect{Proof of Theorem \ref{boundCDFLR} } \la{boundCDFLR_app}

For simplicity of notations, define $F (z, \se) = \Pr \{ \varphi_n \leq z \mid \se \}$ and $G (z, \se) = \Pr \{ \varphi_n \geq z \mid \se \}$.
By the assumption of the theorem, {\small $\f{ f_n( X_1, \cd, X_n; \se) }{f_n( X_1, \cd, X_n; \wh{\se}_n)}  = \Lm (\varphi_n , \wh{\se}_n,
\se)$}.  By virtue of Theorem \ref{THm1}, we have  \bee &  & \Pr \li \{ \f{ f_n( X_1, \cd, X_n; \se) }{f_n( X_1, \cd, X_n; \wh{\se}_n)} \leq
\f{\al}{2}, \; \wh{\se}_n \leq \se, \; \mid \se \ri \} = \Pr \li \{ \Lm (\varphi_n , \wh{\se}_n, \se) \leq
\f{\al}{2}, \; \wh{\se}_n \leq \se \mid \se \ri \}\\
&  & =  \Pr \li \{ \Lm (\varphi_n , \wh{\se}_n, \se)  \leq \f{\al}{2}, \; \wh{\se}_n \leq \se \mid \se \ri \}
\leq  \Pr \li \{  F (\varphi_n, \se)  \leq \f{\al}{2}, \; \wh{\se}_n \leq \se \mid \se \ri \}\\
&  & \leq  \Pr \li \{ F (\varphi_n, \se)  \leq \f{\al}{2} \mid \se \ri \} \leq \f{\al}{2} \eee for any $\se \in \Se$. This proves (\ref{RB1}).
Similarly, for any $\se \in \Se$, \bee &  & \Pr \li \{   \f{ f_n( X_1, \cd, X_n; \se) }{f_n( X_1, \cd, X_n; \wh{\se}_n)}  \leq \f{\al}{2}, \;
\wh{\se}_n \geq \se \mid \se \ri \} = \Pr \li \{ \Lm (\varphi_n , \wh{\se}_n, \se) \leq
\f{\al}{2}, \; \wh{\se}_n \geq \se \mid \se \ri \}\\
&  & =   \Pr \li \{ \Lm (\varphi_n , \wh{\se}_n, \se)  \leq \f{\al}{2}, \; \wh{\se}_n \geq \se \mid \se \ri \} \leq \Pr \li \{ G (\varphi_n,
\se)  \leq \f{\al}{2}, \; \wh{\se}_n \geq \se \mid \se \ri \}\\
&  & \leq  \Pr \li \{ G (\varphi_n, \se)  \leq \f{\al}{2} \mid \se \ri \} \leq \f{\al}{2}, \eee which establishes (\ref{RB2}).  To show
(\ref{RB3}), making use of (\ref{RB1}) and (\ref{RB2}), we have
 \bee &  & \Pr \li \{ \f{ f_n( X_1, \cd, X_n; \se) }{f_n( X_1, \cd, X_n; \wh{\se}_n)} \leq \f{\al}{2} \mid \se \ri \}\\
&  & =  \Pr \li \{  \f{ f_n( X_1, \cd, X_n; \se) }{f_n( X_1, \cd, X_n; \wh{\se}_n)} \leq \f{\al}{2}, \; \wh{\se}_n \leq \se \mid \se \ri \}
 + \Pr \li \{  \f{ f_n( X_1, \cd, X_n; \se) }{f_n( X_1, \cd, X_n; \wh{\se}_n)} \leq \f{\al}{2}, \; \wh{\se}_n \geq \se \mid \se \ri \}\\
&  & \leq  \f{\al}{2} + \f{\al}{2} = \al \eee for any $\se \in \Se$. To show (\ref{LBA}), making use of (\ref{RB1}), we have  that \bee &  & \Pr
\li \{  \f{ \sup_{\vse \in \mscr{S}} f_n( X_1, \cd, X_n; \vse) }{\sup_{\vse \in \Se} f_n( X_1, \cd, X_n; \vse)} \leq \f{\al}{2}, \; \wh{\se}_n
\leq \inf
\mscr{S} \mid \se \ri \}\\
&  & \leq \Pr \li \{  \f{ \sup_{\vse \in \mscr{S}} f_n( X_1, \cd, X_n; \vse) }{\sup_{\vse \in \Se} f_n( X_1, \cd, X_n; \vse)} \leq \f{\al}{2},
\; \wh{\se}_n \leq \se \mid \se \ri \}\\
&  & \leq \Pr \li \{  \f{ f_n( X_1, \cd, X_n; \se) }{\sup_{\vse \in \Se} f_n( X_1, \cd, X_n; \vse)} \leq \f{\al}{2}, \; \wh{\se}_n \leq \se \mid
\se \ri \}\\
&  & = \Pr \li \{  \f{ f_n( X_1, \cd, X_n; \se) }{f_n( X_1, \cd, X_n; \wh{\se}_n)} \leq \f{\al}{2}, \; \wh{\se}_n \leq \se \mid \se \ri \} \leq
\f{\al}{2} \eee for any $\se \in \mscr{S}$.   To show (\ref{LBB}),  making use of (\ref{RB2}), we have that \bee &  & \Pr \li \{  \f{ \sup_{\vse
\in \mscr{S}} f_n( X_1, \cd, X_n; \vse) }{\sup_{\vse \in \Se} f_n( X_1, \cd, X_n; \vse)} \leq \f{\al}{2}, \; \wh{\se}_n \geq \sup
\mscr{S} \mid \se \ri \}\\
&  & \leq \Pr \li \{  \f{ \sup_{\vse \in \mscr{S}} f_n( X_1, \cd, X_n; \vse) }{\sup_{\vse \in \Se} f_n( X_1, \cd, X_n; \vse)} \leq \f{\al}{2},
\; \wh{\se}_n \geq \se \mid \se \ri \}\\
&  & \leq \Pr \li \{  \f{ f_n( X_1, \cd, X_n; \se) }{\sup_{\vse \in \Se} f_n( X_1, \cd, X_n; \vse)} \leq \f{\al}{2}, \; \wh{\se}_n \geq \se \mid
\se \ri \}\\
&  & = \Pr \li \{  \f{ f_n( X_1, \cd, X_n; \se) }{f_n( X_1, \cd, X_n; \wh{\se}_n)} \leq \f{\al}{2}, \; \wh{\se}_n \geq \se \mid \se \ri \} \leq
\f{\al}{2} \eee for any $\se \in \mscr{S}$.  To show (\ref{LBC}), we use (\ref{RB3}) to conclude that \bee \Pr \li \{  \f{ \sup_{\vse \in
\mscr{S}} f_n( X_1, \cd, X_n; \vse) }{\sup_{\vse \in \Se} f_n( X_1, \cd, X_n; \vse)} \leq \f{\al}{2} \mid \se \ri \}  &  \leq & \Pr \li \{  \f{
f_n( X_1, \cd, X_n; \se) }{\sup_{\vse \in \Se} f_n( X_1, \cd,
X_n; \vse)} \leq \f{\al}{2} \mid \se \ri \} \\
&  =  & \Pr \li \{  \f{ f_n( X_1, \cd, X_n; \se) }{f_n( X_1, \cd, X_n; \wh{\se}_n)} \leq \f{\al}{2} \mid \se \ri \} \leq \al  \eee for any $\se
\in \mscr{S}$.  This completes the proof of the theorem.

\section{Proof of Theorem \ref{THm2} } \la{THm2_app}

Note that {\small $\prod_{i = 1}^n f_X (x_i, \se) =  \li [ \prod_{i
= 1}^n h (x_i) \ri ] \times \exp \li ( \eta (\se) \sum_{i = 1}^n T (
x_i ) - n A (\se) \ri )$}.  By the assumption that $\f{d
\eta(\se)}{d \se}$ is positive for $\se \in \Se$, we have that the
likelihood ratio
\[
\Lambda (z, \se_0, \se_1) = \li [ \f{  \exp \li ( \eta (\se_1) z - A
(\se_1)  \ri )  } {  \exp \li ( \eta (\se_0) z  - A (\se_0) \ri ) }
\ri ]^n
\]
is an increasing function of $z \in \Se$ provided that $\se_0 < \se_1$. Applying Theorem \ref{THm1} with $\vse (z) = z$, we have
\[
\Pr \{ \wh{\bs{\se}} \geq z \mid \se \} \leq \li [ \f{  \exp \li (
\eta (\se) z - A (\se)  \ri )  } {  \exp \li ( \eta (z)  z  -  A (z)
\ri ) } \ri ]^n \times \Pr \{ \wh{\bs{\se}} \geq z \mid z \} =
\mscr{M} (z, \se) \times \Pr \{ \wh{\bs{\se}} \geq z \mid z \}
\]
for $z \in \Se$ no less than $\se \in \Se$. Similarly, $\Pr \{ \wh{\bs{\se}} \leq z \mid \se \} \leq \mscr{M} (z, \se) \times \Pr \{
\wh{\bs{\se}} \leq z \mid z \}$ for $z \in \Se$ no greater than $\se \in \Se$.  It remains to show statements (i)--(v) under the additional
assumption that $\f{ A^\prime ( \se ) }{ \eta^\prime ( \se ) } = \se$. For simplicity of notations, define $w(z, \se) = \exp \li ( \eta (\se) z
- A (\se)  \ri )$.  Since $\f{d \eta (\se)}{d \se} > 0$ and $\f{ A^\prime ( \se ) }{ \eta^\prime ( \se ) } = \se$ for $\se \in \Se$, we have
that
\[
\f{ d w(z, \se)  } { d \se } =  (z - \se) w(z, \se)  \f{d \eta
(\se)}{d \se},
\]
which is positive for $\se < z$ and negative for $\se > z$.  This
implies that $w (z, \se)$ is monotonically increasing with respect
to $\se$ less than $z$ and monotonically decreasing with respect to
$\se$ greater than $z$. Therefore, $\wh{\bs{\se}}$ must be a
maximum-likelihood estimator of $\se$.

Let $\psi (.)$ be the inverse function of $\eta(.)$ such that \be
\la{def89} \eta (\psi (\ze) ) = \ze \ee for $\ze \in \{ \eta (\se):
\se \in \Se \}$.  Define compound function $B(.)$ such that $B
(\zeta) = A (\psi (\zeta))$ for $\ze \in \{ \eta (\se): \se \in \Se
\}$. For simplicity of notations, we abbreviate $\psi (\zeta)$ as
$\psi$ when this can be done without causing confusion. By the
assumption that $\f{ d A ( \se ) }{ d \se } = \se \f{ d \eta ( \se )
}{d \se}$, we have \be \la{def836}
 \f{ \f{ d A (\psi) }{d \psi} } {
\f{ d \eta (\psi) }{d \psi} } = \psi. \ee Using (\ref{def89}),
(\ref{def836}) and the chain rule of differentiation, we have \be
\la{use8} \f{ d B (\zeta ) }{d \zeta} = \f{ d A (\psi) }{d \psi} \f{
d \psi }{ d \zeta} = \f{ \f{ d A (\psi) }{d \psi} } { \f{ d \eta
(\psi) }{d \psi} } \f{ d \eta (\psi) }{d \psi} \f{ d \psi }{ d
\zeta} = \f{ \f{ d A (\psi) }{d \psi} } { \f{ d \eta (\psi) }{d
\psi} } \f{ d \eta (\psi) }{d \ze} = \psi \f{ d \ze }{d \ze} = \psi
(\ze). \ee  Putting $\zeta = \eta (\se)$, we have {\small \bee  & &
\bb{E} \li [ \exp \li ( n t ( \wh{\bs{\se}}) \ri ) \ri ]  = \bb{E}
\li [ \exp \li ( t \sum_{i = 1}^n T ( X_i ) \ri )  \ri ] =  \int \cd
\int \prod_{i = 1}^n \li [ h (x_i)  \exp \li ( ( \zeta +
t ) T ( x_i ) -  B (\zeta) \ri ) \ri ] dx_1 \cd dx_n\\
& = &  \exp \li (  n B (\zeta + t) - n B (\zeta) \ri ) \int \cd \int
\prod_{i = 1}^n \li [ h (x_i)
\exp \li ( ( \zeta + t ) T ( x_i ) - B (\zeta + t) \ri ) \ri ] dx_1 \cd dx_n\\
& = &  \exp \li (  n B (\zeta + t) - n B (\zeta) \ri ). \eee} By
virtue of (\ref{use8}), the derivative of $n B (\zeta + t) - n B
(\zeta)$ with respect to $t$ is
\[
n \f{d B (\zeta + t) }{d t} = n \psi ( \zeta + t),
\]
which is equal to $n \psi ( \zeta ) = n \se$ for $t = 0$.  Thus,
$\bb{E} [ \wh{\bs{\se}} ] = \se$, which implies that $\wh{\bs{\se}}$
is also an unbiased estimator of $\se$.  This proves statement (i).

Again by virtue of (\ref{use8}), the derivative of $- t n z + n B
(\zeta + t) - n B (\zeta)$ with respect to $t$ is
\[
-  n z +  n \f{d B (\zeta + t) }{d t} = - n z + n \psi ( \zeta + t),
\]
which is equal to $0$ for $t$ such that $\psi ( \zeta + t ) = z$ or
equivalently, $\zeta + t = \eta (z)$, which implies $t = \eta (z) -
\eta (\se )$.  Since $\bb{E} \li [ \exp \li ( n t ( \wh{\bs{\se}} -
z) \ri ) \ri ]$ is a convex function of $t$, its infimum with
respect to $t \in \bb{R}$ is attained at $t = \eta (z) - \eta (\se
)$. It follows that \bee & & \inf_{t \in \bb{R} } \bb{E} \li [ \exp
\li ( n t ( \wh{\bs{\se}} - z) \ri ) \ri ] = \inf_{t \in
\bb{R}} \exp \li ( - t n z +  n B (\zeta + t) - n B (\zeta) \ri )\\
&  =  & \exp \li ( - [ \eta (z) - \eta (\se ) ]  n z +  n B ( \eta
(z) ) - n B (\zeta) \ri )  =  \exp \li ( - [ \eta (z) - \eta (\se )
]  n z +  n A ( z) - n
A (\se) \ri )\\
& = &  \li [ \f{ \exp \li ( \eta (\se )   z -  A (\se) \ri ) } {
\exp \li ( \eta (z ) z - A (z) \ri ) } \ri ]^n = \mscr{M} (z, \se ).
\eee

Now, consider the monotonicity of $\mscr{M} (z, \se)$ with respect
to $\se$ as described by statement (iii) of the theorem.  To show
$\mscr{M} (z, \se_2) \leq \mscr{M} (z, \se_1)$ for $\se_2
> \se_1 \geq z$, note that
\[
\mscr{M} (z, \se_2) = \li [ \f{ w (z, \se_2) } { w (z, z) } \ri ]^n
\leq \li [ \f{ w (z, \se_1) } { w (z, z) } \ri ]^n = \mscr{M} (z,
\se_1),
\]
where the inequality is due to the fact that $w(z, \se)$ is
non-increasing with respect to $\se$ no less than $z$.  On the other
hand,  to show $\mscr{M} (z, \se_1) \leq \mscr{M} (z, \se_2)$ for
$\se_1 < \se_2 \leq z$,  note that
\[
\mscr{M} (z, \se_1) = \li [ \f{ w (z, \se_1) } { w (z, z) } \ri ]^n
\leq \li [ \f{ w (z, \se_2) } { w (z, z) } \ri ]^n = \mscr{M} (z,
\se_2),
\] where the inequality is due to  the fact that $w(z, \se)$ is
non-decreasing with respect to $\se$ no greater than $z$. This
justifies statement (iii) of the theorem.

Next, consider the monotonicity of $\mscr{M} (z, \se)$ with respect
$z$ as described by statement (iv) of the theorem.  To show
$\mscr{M} (z_2, \se) \leq \mscr{M} (z_1, \se)$ for $z_2
> z_1 \geq \se$, it is sufficient to note that
\[
\mscr{M} (z_2, \se) = \li [ \f{ w (z_2, \se) } { w (z_2, z_2) } \ri
]^n \leq \li [ \f{ w (z_2, \se) } { w (z_2, z_1) } \ri ]^n \leq \li
[ \f{ w (z_1, \se) } { w (z_1, z_1) } \ri ]^n = \mscr{M} (z_1, \se),
\]
where the first inequality is due to the fact that $w(z, \se)$ is
non-decreasing with respect to $\se$ no greater than $z$ and the
second one is due to the assumption that the likelihood ratio
$\Lambda (z, \se_0, \se_1) = \li [ \f{ w(z, \se_1) }{ w(z, \se_0) }
\ri ]^n$ is non-decreasing with respect to $z$.  On the other side,
to show $\mscr{M} (z_2, \se) \geq \mscr{M} (z_1, \se)$ for $z_1 <
z_2 \leq \se$, it suffices to observe that
\[
\mscr{M} (z_1, \se) = \li [ \f{ w (z_1, \se) } { w (z_1, z_1) } \ri
]^n \leq \li [ \f{ w (z_1, \se) } { w (z_1, z_2) } \ri ]^n \leq \li
[ \f{ w (z_2, \se) } { w (z_2, z_2) } \ri ]^n = \mscr{M} (z_2, \se),
\]
where the first inequality is due to the fact that $w(z, \se)$ is
non-increasing with respect to $\se$ no less than $z$ and the second
one is due to the assumption that the likelihood ratio $\Lambda (z,
\se_0, \se_1) = \li [ \f{ w(z, \se_1) }{ w(z, \se_0) } \ri ]^n$ is
non-decreasing with respect to $z$. Statement (iv) of the theorem is
thus proved.

Finally, in order to show statement (v), notice that, in the course
of proving that $\wh{\bs{\se}}$ is an unbiased estimator of $\se$,
we have shown that $\bb{E} [ T(X) - \se ]  = 0$.  Hence, applying
the Berry-Essen inequality and Lyapounov's inequality, we have that
both (\ref{ineq8}) and (\ref{ineq88}) are true.

\sect{Proof of Theorem \ref{LRB_CI} } \la{LRB_CI_app}

For simplicity of notations, define $F_{\bs{\varphi}} (z, \se) = \Pr
\{ \bs{\varphi} \leq z \mid \se \}$ and $G_{\bs{\varphi}} (z, \se) =
\Pr \{ \bs{\varphi} \geq z \mid \se \}$.  By the assumption of the
theorem, we have \be \la{cr1} F_{\bs{\varphi}} (z, \se) \leq \bb{F}
(z, \se) \qu \tx{for} \; z \leq \se, \ee \be \la{cr2}
G_{\bs{\varphi}} (z, \se) \leq \bb{G} (z, \se)  \qu \tx{for} \; z
\geq \se. \ee Making use of (\ref{cr1}), the assumption that $\bb{F}
(z, \se)$ is non-increasing with respect to $\se \geq z$,  and the
assumption that $\{ \bb{F} ( \bs{\varphi}, U (\bs{\varphi}, \de) )
\leq \f{\de}{2}, \; \bs{\varphi} \leq U (\bs{\varphi}, \de) \}$ is a
sure event, we have {\small \bee &
 & \{ U (\bs{\varphi}, \de) \leq \se \} = \{ \bs{\varphi} \leq U
(\bs{\varphi}, \de) \leq \se, \; \bb{F} (\bs{\varphi}, U
(\bs{\varphi}, \de) )  \leq \f{\de}{2} \}\\
&  & \subseteq \{ \bs{\varphi} \leq U (\bs{\varphi}, \de) \leq \se,
\; \bb{F} (\bs{\varphi}, \se ) \leq
\f{\de}{2} \}\\
 &  &  \subseteq \{ \bs{\varphi} \leq U (\bs{\varphi}, \de) \leq \se,
\; F_{\bs{\varphi}} (\bs{\varphi}, \se ) \leq \f{\de}{2} \}
\subseteq \{ F_{\bs{\varphi}} (\bs{\varphi}, \se ) \leq \f{\de}{2}
\}, \eee} which implies that $\Pr \{ U (\bs{\varphi}, \de) \leq \se
\} \leq \Pr \{ F_{\bs{\varphi}} (\bs{\varphi}, \se ) \leq \f{\de}{2}
\} \leq \f{\de}{2}$.
  On the other hand,  Making use of
(\ref{cr2}), the assumption that $\bb{G} (z, \se)$ is non-decreasing
with respect to $\se \leq z$,  and the assumption that $\{ \bb{G} (
\bs{\varphi}, L (\bs{\varphi}, \de) )  \leq \f{\de}{2}, \;
\bs{\varphi} \geq L (\bs{\varphi}, \de) \}$ is a sure event, we have
{\small \bee &
 & \{ L (\bs{\varphi}, \de) \geq \se \} = \{ \bs{\varphi} \geq L
(\bs{\varphi}, \de) \geq \se, \; \bb{G} (\bs{\varphi}, L
(\bs{\varphi}, \de) ) \leq \f{\de}{2} \}\\
&  & \subseteq \{ \bs{\varphi} \geq L (\bs{\varphi}, \de) \geq \se,
\; \bb{G} (\bs{\varphi}, \se ) \leq
\f{\de}{2} \}\\
 &  &  \subseteq \{ \bs{\varphi} \geq L (\bs{\varphi}, \de) \geq \se,
\; G_{\bs{\varphi}} (\bs{\varphi}, \se ) \leq \f{\de}{2} \}
\subseteq \{ G_{\bs{\varphi}} (\bs{\varphi}, \se ) \leq \f{\de}{2}
\}, \eee} which implies that $\Pr \{ L (\bs{\varphi}, \de) \geq \se
\} \leq \Pr \{ G_{\bs{\varphi}} (\bs{\varphi}, \se ) \leq \f{\de}{2}
\} \leq \f{\de}{2}$.  Finally, by virtue of the established fact
that $\Pr \{ U (\bs{\varphi}, \de) \leq \se \} \leq \f{\de}{2}$ and
$\Pr \{ L (\bs{\varphi}, \de) \geq \se \} \leq \f{\de}{2}$, we have
$\Pr \{ L (\bs{\varphi}, \de) < \se < U (\bs{\varphi}, \de) \mid \se
\} \geq 1 - \Pr \{ U (\bs{\varphi}, \de) \leq \se \} - \Pr \{ L
(\bs{\varphi}, \de) \geq \se \} \geq 1 - \f{\de}{2} - \f{\de}{2} = 1
- \de$. This completes the proof of the theorem.

\sect{Proof of Theorem \ref{AST_Bound} } \la{AST_Bound_App}

Let $N_{\mrm{b}}(\de)$ be the minimum sample size to ensure that
\[
\Pr \{ \ovl{X}_n \geq \mu + \vep \} \leq \f{\de}{2}, \qqu \Pr \{
\ovl{X}_n \leq \mu - \vep \} \leq \f{\de}{2}.
\]
Since $\Pr \{ | \ovl{X}_n - \mu | \geq \vep \}$ equals the summation
of $\Pr \{ \ovl{X}_n \geq \mu + \vep \}$ and $\Pr \{ \ovl{X}_n \leq
\mu - \vep \}$,  we have that $\Pr \{ | \ovl{X}_n - \mu | \geq \vep
\} \leq \de$ implies $\Pr \{ \ovl{X}_n \geq \mu + \vep \} \leq \de$
and $\Pr \{ \ovl{X}_n \leq \mu - \vep \} \leq \de$. Consequently,
\[
N_{\mrm{a}} (\de)
> N_{\mrm{b}}(2 \de).
\]
Since $\Pr \{ \ovl{X}_n \geq \mu + \vep \} \leq \f{\de}{2}$ and $\Pr
\{ \ovl{X}_n \leq \mu - \vep \} \leq \f{\de}{2}$ together imply $\Pr
\{ | \ovl{X}_n - \mu | \geq \vep \} \leq \de$, we have \[
N_{\mrm{a}} (\de) < N_{\mrm{b}}(\de). \] Therefore, $N_{\mrm{b}}(2
\de) < N_{\mrm{a}} (\de) < N_{\mrm{b}}(\de)$. We claim that
$\lim_{\de \to 0} \f{ N_{\mrm{c}}(\de)  } { N_{\mrm{b}}(\de) } = 1$.
To show this claim, we define {\small \[ Q^+ = \f{\ln \Pr \{
\ovl{X}_n \geq \mu + \vep \}}{n} \]} and {\small \[ Q^- = \f{\ln \Pr
\{ \ovl{X}_n \leq \mu - \vep \}}{n}.  \]}  Then, \[ Q^+ < 0, \qu Q^-
< 0, \qu  \ln \mcal{F}(\mu - \vep, \mu) < 0, \qu \ln \mcal{G}(\mu +
\vep, \mu) < 0 \]
 and {\small
 \[
 N_{\mrm{b}}(\de) = \max \li \{ \f{
\ln \f{\de}{2} } { Q^+ }, \; \f{ \ln \f{\de}{2} } { Q^-} \ri \}. \]}
It follows that {\small \bee &  & \lim_{\de \to 0} \f{ N_{\mrm{c}}(
\de) } { N_{\mrm{b}}(\de) } = \lim_{\de \to 0} \f{1} {
N_{\mrm{b}}(\de) } \times \max \li \{ \f{ \ln \f{\de}{2} } { \ln
\mcal{F}(\mu - \vep, \mu) }, \f{ \ln \f{\de}{2} } { \ln \mcal{G}(\mu
+ \vep, \mu) } \ri \}\\
&  &  = \lim_{\de \to 0} \f { \max \{ Q^+, Q^- \} } { \max \{ \ln
\mcal{F}(\mu - \vep, \mu), \ln \mcal{G}(\mu + \vep, \mu) \} }. \eee}
By Chernoff's theorem,
\[
\lim_{\de \to 0} Q^+ = \ln \mcal{G}(\mu + \vep, \mu), \qqu \lim_{\de
\to 0} Q^- = \ln \mcal{F}(\mu - \vep, \mu) \]
 and consequently,
\[
\lim_{\de \to 0} \max \{ Q^+, Q^- \} = \max \{ \ln \mcal{F}(\mu -
\vep, \mu), \; \ln \mcal{G}(\mu + \vep, \mu) \} \] and the claim
follows. Using the established claim, we have {\small \bee \lim_{\de
\to 0} \f{ N_{\mrm{b}}(\de) } { N_{\mrm{b}}(2 \de) }  =  \lim_{\de
\to 0} \li [ \f{ N_{\mrm{b}}(\de) } { N_{\mrm{c}}( \de) } \times \f{
\max \li \{ \f{ \ln \f{\de}{2} } { \ln \mcal{F}(\mu + \vep, \mu) },
\f{ \ln \f{\de}{2} } { \ln \mcal{G}(\mu - \vep, \mu) } \ri \} }  {
\max \li \{ \f{ \ln \de } { \ln \mcal{F}(\mu + \vep, \mu) }, \f{ \ln
\de } { \ln \mcal{G}(\mu - \vep, \mu) } \ri \} }  \times  \f{
N_{\mrm{c}}( 2 \de)  } { N_{\mrm{b}}(2 \de) } \ri ] = \lim_{\de \to
0} \f{ \ln \f{\de}{2} } { \ln \de } = 1. \eee} Recalling
$N_{\mrm{b}}(2 \de) < N_{\mrm{a}} (\de) < N_{\mrm{b}}(\de)$,  we can
conclude that $\lim_{\de \to 0} \f{ N_{\mrm{b}}(\de)  } {
N_{\mrm{a}} (\de) } = 1$.  Finally, recalling the established claim
that $\lim_{\de \to 0} \f{ N_{\mrm{c}}(\de)  } { N_{\mrm{b}} (\de) }
= 1$, the proof of the theorem is thus completed.

\end{document}